\documentclass[a4paper,11pt]{article} %
\usepackage[T1]{fontenc} % Better handling of accented characters
\usepackage[utf8]{inputenc} % Allows UTF-8 characters
\usepackage[portuguese,english]{babel} % English is the default 
\usepackage[a4paper,margin=1in]{geometry} % Adjust page layout
\usepackage{amsmath,amsthm,amssymb,amsfonts}
\usepackage{mathtools} % loads the amsmath package automatically.
\usepackage{bm}% used to get bold symbos; command *\bm{xxxx}*
\usepackage{mathrsfs} % try the command *$\mathscr{H}$*
\usepackage{changes} % Used to make changes in a text, e.g., with the command *\deleted{xxxx}*
\usepackage{url} % used to write *email adresses*, *hypertexts links*, etc.
\usepackage[toc,page]{appendix} %-> a opção ``page'' permite criar um ``environment'' e a opção ``toc'' 
\usepackage{enumerate}%The package adds an optional argument to the enumerate environment which determines the style in which the counter is printed
\usepackage{enumitem}% To customize the label in itemize
\usepackage{cases}% To use separate labels in the "cases" enviorment
\usepackage{enumitem}
\usepackage{tcolorbox}
\usepackage{mdframed}
\usepackage{graphicx} %\usepackage[dvips]{graphicx}  %\usepackage[pdftex]{graphicx}
\usepackage{longtable} %Allow tables to flow over page boundaries
\usepackage[all]{xy} % para usar xy-pic, fazer diagramas.
\usepackage{newfloat} %The package offers the command \DeclareFloatingEnvironment
\DeclareFloatingEnvironment[fileext=frm,placement={ht},name=Frame]{algfloat}
\usepackage{caption}%provides many ways to customise the captions in floating environments like figure and table
\usepackage{placeins}%Defines a \FloatBarrier command, beyond which floats may not pass; useful, for example, to ensure all %floats for a section appear before the next \section command
\usepackage{tikz}%TikZ is probably the most complex and powerful tool to create graphic elements in LaTeX.
\usepackage{subfig}%The package provides support for the manipulation and reference of small or ‘sub’ figures and tables within a single figure or table environment
\usepackage{pgfplots}%PGFPlots draws high-quality function plots
\usepackage{subcaption} 

\usepackage{xcolor} %\usepackage[usenames,dvipsnames,svgnames,table]{xcolor}%The package starts from the basic facilities of the color package, and provides easy driver-independent access to several kinds of color tints
\usepackage{array}%An extended implementation of the array and tabular environments which extends the options for column formats
\usepackage{multirow} % Introd. paper com João e Benar
\usepackage{booktabs}%The package enhances the quality of tables in LaTeX, providing extra commands as well as behind-the-scenes optimisation
\usepackage{grffile}%The original package extended the file name processing of package graphics to support a larger range of file names. 
\newtheorem{theorem}{Theorem}[section] %[section] restarts the theorem counter at every new section
\newtheorem{lemma}[theorem]{Lemma}
\newtheorem{corollary}[theorem]{Corollary}
\newtheorem{proposition}[theorem]{Proposition}
\newtheorem{remark}[theorem]{Remark}
\newtheorem{definition}[theorem]{Definition}
\newtheorem{assumption}[theorem]{Assumption}

\newtheorem{example}[theorem]{Example}

\newcommand{\Z}{\mathbb{Z}}

\newcommand{\R}{\mathbb{R}}

\newcommand{\ER}{(-\infty, \infty]}%{\overline{\mathbb{R}}}
\newcommand{\inner}[2]{\langle #1,#2\rangle} % [2] number of arguments. Modo de uso: \inner{a}{b}
\newcommand{\Inner}[2]{\left\langle #1,#2\right\rangle}
\newcommand{\norm}[1]{\|{#1}\|}
\newcommand{\normq}[1]{\|{#1}\|^2}

\newcommand{\veps}{{\varepsilon}}

\newcommand{\what}[1]{\widehat #1}

\newcommand{\HH}{\R^n} % Espaço de Hilbert

\newcommand{\comenta}[1]{} % comentar textos - no texworks pode ser feito com CTRL + SHFIT + ]
\newcommand{\set}[1]{\{#1\}}

\newcommand{\seq}[1]{\left(#1\right)}

\newcommand{\lab}[1]{\label{#1}}

\newcommand{\bprop}{\begin{proposition}}
\newcommand{\eprop}{\end{proposition}}
\newcommand{\blemm}{\begin{lemma}}
\newcommand{\elemm}{\end{lemma}}
\newcommand{\bdefi}{\begin{definition}}
\newcommand{\edefi}{\end{definition}}
\newcommand{\btheo}{\begin{theorem}}
\newcommand{\etheo}{\end{theorem}}
\newcommand{\bproo}{\begin{proof}}
\newcommand{\eproo}{\end{proof}}
\newcommand{\brema}{\begin{remark}}
\newcommand{\erema}{\end{remark}}
\newcommand{\bitem}{\begin{itemize}}
\newcommand{\eitem}{\end{itemize}}
\newcommand{\bexam}{\begin{example}}
\newcommand{\eexam}{\end{example}}
\newcommand{\bassu}{\begin{assumption}}
\newcommand{\eassu}{\end{assumption}}
\newcommand{\bcoro}{\begin{corollary}}
\newcommand{\ecoro}{\end{corollary}}

\newcommand{\benum}{\begin{enumerate}[label = \emph{(\alph*)}]}
\newcommand{\eenum}{\end{enumerate}}

\newcommand{\mgap}{\vspace{.1in}}

\newcommand{\beq}{\begin{equation}}
\newcommand{\eeq}{\end{equation}}

\DeclareMathOperator{\Dom}{Dom} % \newcommand is in general used for much more complex situations

\DeclareMathOperator{\dist}{dist}
\DeclareMathOperator{\Gap}{Gap}
\usepackage[ruled]{algorithm2e}%linesnumbered, 
\SetKwInput{KwInput}{Input}
\SetKwInput{KwOutput}{Global Variables for Function}
\definechangesauthor[name={M. Marques Alves}, color=blue]{M}
\definechangesauthor[name={Kangming Chen}, color=red]{C}
\definechangesauthor[name={Ellen H. Fukuda}, color=green]{E}
\definecolor{myred}{RGB}{255,50,50}         % revision version
\definecolor{myblack}{RGB}{0,0,0}           % normal version

\begin{document}

\title{An inertial iteratively regularized extragradient method for bilevel variational inequality problems}
\author{ M. Marques Alves%
\thanks{
Departamento de Matem\'atica,
Universidade Federal de Santa Catarina,
Florian\'opolis, 88040-900, Brazil ({\tt maicon.alves@ufsc.br}).
The work of this author was partially supported by CNPq, Grant 308036/2021-2, and Fundação de Amparo a
Pesquisa e Inovação do Estado de Santa Catarina (FAPESC), Edital 21/2024, Grant 2024TR002238.
}
\and
Kangming Chen%
\thanks{
Graduate School of Informatics, Kyoto University, Sakyo-ku Yoshida-Honmachi, Kyoto
606–8501, Japan ({\tt kangming@amp.i.kyoto-u.ac.jp}).}
\and
Ellen H. Fukuda%
\thanks{
Graduate School of Informatics, Kyoto University, Sakyo-ku Yoshida-Honmachi, Kyoto
606–8501, Japan ({\tt ellen@i.kyoto-u.ac.jp}). This author is supported by 
Japan Society for the Promotion of Science, Grant-in-Aid for Scientific Research (C) (JP25K15002) and Grant-in-Aid for Scientific Research (B) (JP25K03082).}
}

%\date{May 9, 2014}
%\date{}

\maketitle
%\tableofcontents

\begin{abstract}
We study a bilevel variational inequality problem where the feasible set is itself the solution set of another variational inequality. Motivated by the difficulty of computing projections onto such sets, we consider a regularized extragradient method, as proposed by Samadi and Yousefian (2025), which operates over a simpler constraint set. 
Building on this framework, we introduce an inertial variant (called IneIREG) that incorporates momentum through extrapolation steps. We establish iteration-complexity bounds for the general (non-strongly monotone) case under both constant and diminishing regularization, and derive improved results under strong monotonicity assumptions. Our analysis extends and refines the results of the previous work by capturing both inertial and regularization effects within a unified framework. Preliminary numerical experiments are also presented to illustrate the behavior of the proposed method.
\\
\\ 
  2000 Mathematics Subject Classification: 90C33, 49J40, 65K15, 68Q25.
%
%90C25 (Convex Programming)
%90C06 (Large Scale Problems in Mathematical Programming)
%
%47H05 (Monotone Operators and Generalizations).
%47N10 (Applications of operator theory in optimization, convex analysis, mathematical programming, economics)
%47J25  (Iterative procedures involving nonlinear operators)
%    
%49J40  (Variational inequalities)
%49J52 (Nonsmooth Analysis)
%49M15 (Newton-type methods)
%
%65K15  (Numerical methods for variational inequalities and related problems)
%68Q25 (Analysis of algorithms and problem complexity)
%
%90C30 (Nonlinear Programming)
%90C33 (Complementarity and equilibrium problems and variational inequalities)
\end{abstract}

\pagestyle{plain}

\section{Introduction}

We consider the (bilevel) variational inequality problem $ \text{VIP}(H, Q) $:
\begin{align} \lab{eq:prob}
\mbox{Find} \enspace x \in Q \enspace \mbox{such that} \enspace
\Inner{H(x)}{z - x} \geq 0 \enspace \text{for all} \enspace z \in Q,
\end{align}
where $ Q $ is the solution set of the variational inequality $ \text{VIP}(F, X) $, i.e.,
\begin{align} \lab{eq:def.Q}
 Q \coloneqq \set{x \in X \mid \Inner{F(x)}{y - x} \geq 0 \enspace \text{for all} \enspace y \in X}.
\end{align}
(see Section \ref{sec:main} for the precise statement.)
Problem \eqref{eq:prob} arises in different contexts in applications, ranging from (bilevel) optimization to game theory and operations research (see, e.g., \cite[Subsection 1.1]{SamYou25}). 

One of the most widely used algorithms for solving variational inequality problems is the celebrated extragradient (EG) method proposed by Korpelevich~\cite{Kor76}, which, when applied to problem~\eqref{eq:prob}, takes the following form:
\begin{align} \lab{eq:kor.intr}
y_k = P_Q(x_k - \lambda_k H(x_k)), \quad x_{k+1} = P_Q(x_k - \lambda_k H(y_k)),
\end{align}
where $ x_k $ and $ x_{k+1} $ denote the current and next iterates, respectively, $ \lambda_k > 0 $ is the stepsize parameter
and $ P_Q(\cdot) $ denotes the orthogonal projection onto $ Q $. 
The main challenge in implementing \eqref{eq:kor.intr} lies in the fact that the projection 
$ P_Q(\cdot) $ is difficult to compute, since the feasible set
$ Q $, defined in \eqref{eq:def.Q}, is itself the solution set of another variational inequality. A common alternative found in the literature (see, e.g., \cite{SamYou25}) is to adopt a regularization approach, where the operators 
$ F(\cdot) $ and $ H(\cdot) $ 
are combined using a regularization parameter, say, $ \eta_k > 0$.
At each iteration $ k $, one then considers the regularized operator $ F(\cdot) + \eta_k H(\cdot) $,
and applies iteration \eqref{eq:kor.intr} with respect to the feasible set $ X $, which is assumed to be simple enough so that the projection $ P_X (\cdot) $ is computationally tractable.
More precisely, we have
\begin{align} \lab{eq:regkor.intr}
y_k = P_X\Big(x_k - \lambda_k \big(F(x_k)+\eta_k H(x_k)\big)\Big), \quad 
x_{k+1} = P_X\Big(x_k - \lambda_k \big(F(y_k)+\eta_k H(y_k)\big)\Big).
\end{align}
In \cite{SamYou25}, conditions on $ \lambda_k > 0 $ and $ \eta_k > 0 $ are established to obtain iteration-complexity bounds for \eqref{eq:regkor.intr} with respect to the dual gap function; see Definition \ref{def:dualgap} below.

\mgap

In this paper, we propose and analyze the iteration-complexity of an inertial variant 
of \eqref{eq:regkor.intr}\,---see Algorithm \ref{alg:main}.
Proximal algorithms incorporating inertial effects for monotone inclusions and related problems were first introduced in the seminal work~\cite{AlvAtt01}, and have since been further developed in several directions by various authors (see, e.g., \cite{AlvEckGerMel20,AlvMar20,Att21,AttCab20, AttCabChbRia18,BotSedVuo23,CorPey25}, and references therein).
The central idea is to introduce an inertial effect at the current iterate 
$ x_k $ by extrapolating in the direction of the previous displacement $ x_k - x_{k-1} $, leading to the extrapolated point
\begin{align} \lab{eq:alexchorin}
 w_k = x_k + \alpha_k(x_k - x_{k-1}),
\end{align}
where $ \alpha_k \geq 0 $. The subsequent update is then performed using the extrapolated point 
$ w_k $ rather than $ x_k $; see Eqs. \eqref{eq:def.wk} and \eqref{eq:yx} below.

\mgap

Our main algorithm, which can be regarded as an inertial variant of \eqref{eq:regkor.intr}, achieves essentially the same 
iteration-complexity bounds as those established in~\cite{SamYou25}, 
while contributing new insights in the strongly monotone case; see Remark \ref{rem:rockafellar} below. 
In particular, this implies that our theoretical contributions remain relevant even for the non-inertial version of the algorithm, which corresponds to setting  $ \alpha_k \equiv 0 $ in \eqref{eq:alexchorin}.

\mgap

\noindent
{\bf Main contributions.} Our main results can be summarized as follows:
\begin{itemize}
\item[(i)] An inertial iteratively regularized EG (IneIREG) method, namely Algorithm \ref{alg:main} below, for solving the bilevel 
VIP \eqref{eq:prob}. 
\item[(ii)] Propositions \ref{pro:optm01} and \ref{pro:feas}: the first results on the iteration-complexity of Algorithm \ref{alg:main}. In Theorems \ref{the:main} and \ref{the:ofcons02} (see also \ref{the:ofcons}), we prove specialized results for 
diminishing and constant regularization parameter $ \eta_k > 0 $, respectively. 
\item[(iii)] Improved results for the strongly monotone case; we refer the reader to Remarks \ref{rem:osher}, \ref{rem:rockafellar}, \ref{rem:shor} and \ref{rem:suasuna} for additional discussions.
\end{itemize} 

\mgap

\noindent
{\bf Related works.} Numerical schemes for bilevel VIPs and related applications in optimization and game theory have been the subject of attention of several recent works, including~\cite{FacPanScuLam14,HieMou21,KauYou21,SamYou25,ThoTriLiDon20}. 
We refer the reader to~\cite{SamYou25} for a comprehensive overview of the current state of research on bilevel VIPs and bilevel optimization, including applications to
\emph{optimal solution selection problems}, 
\emph{optimal Nash equilibrium seeking problems} and 
\emph{generalized Nash equilibrium problems} (GNEPs).

\mgap

\noindent
{\bf General notation and basic definitions.}
Throughout this paper, $\HH$ denotes the Euclidean $n$-dimensional vector space
with inner product $\inner{\cdot}{\cdot}$ and induced norm $\norm{\cdot} = \sqrt{\inner{\cdot}{\cdot}}$. The transpose of a matrix $A$ is given by $A^\top$.
For a nonempty closed and convex subset $ X $ of $ \HH $, the \emph{orthogonal projection} onto $ X $ is denoted by 
$ P_X(\cdot) $, i.e., for $ x \in \HH $, the vector $ P_X(x) $ is the unique element in $ X $ such that 
$\norm{x - P_X(x)}\leq \norm{x - y}$ for all $ y\in X $. 
The distance of $x \in \HH$ to the set $X$ is denoted by $\dist(x, X)$.
The set $ X $ is said to be \emph{simple} if 
$ P_X(x) $ is easy to compute for all $ x \in \HH $.
A single-valued map $ F \colon \Dom F \to \HH $, where the domain $ \Dom F $ is a subset of $ \HH $, is said to be a 
\emph{monotone operator} if $ \inner{F(x) - F(y)}{x - y} \geq 0 $ for all $ x, y \in \Dom F $. Moreover, $ F(\cdot) $ is said to
be $ L_F $-\emph{Lipschitz continuous} if $ L_F > 0 $ and 
$ \norm{F(x) - F(y)} \leq L_F\norm{x - y} $ for all $ x, y \in \Dom F$.

By $ D_X  > 0 $, we denote  the \emph{diameter} of the compact and convex set $ X \subset \HH $, i.e.,
$ D_X \coloneqq \sup_{x,y\in X}\,\norm{x - y} $.
Define also
\begin{align}\lab{eq:def.ch}
C_H\coloneqq \sup_{x\in X}\,\norm{H(x)} \quad \mbox{and} \quad B_H \coloneqq \sup_{x\in Q}\,\norm{H(x)},
\end{align}
where $ Q $ is as in \eqref{eq:def.Q}.
By $ \Z_+ $ and $ \R_+ $, we denote the set of nonnegative integers and (nonnegative) reals, respectively. We also define the following function, which is the equivalent of a merit function for VI problems.

\bdefi[Dual gap function] \lab{def:dualgap}
The \emph{dual gap function} $\emph{Gap}(\cdot, H, Q) \colon \HH \to \ER $ associated to \emph{VIP}(H, Q) is defined, at every $ z \in \HH $, as
\[
\emph{Gap}(z, H, Q) \coloneqq \sup_{x\in Q}\,\Inner{H(x)}{z - x}.
\]
In the same way, we also define the \emph{dual gap function} $\emph{Gap}(\cdot, F, X)$ 
associated to \emph{VIP}(F, X).
\edefi

\mgap

Note that $\Gap(z, H, Q) \geq 0$ when $z \in Q$, and $\Gap(z, H, Q) = 0$ if and only if $z$ is a solution of VIP($H, Q$).
The proof of the following result follows from~\cite[Theorem 2]{SamYou25}. For the sake of completeness, we present a proof here.

\bprop \lab{pro:optm02}
For all $ y \in \HH $,
\begin{align} \lab{eq:bh}
 \emph{Gap}(y, H, Q) \geq -B_H \dist(y, Q).
\end{align}
Moreover, if $ Q $ is $ \sigma $-weakly sharp
\footnote{there exist $ \sigma > 0 $ and $ \mathcal{M}\geq 1 $ such that 
$ \inner{F(x)}{y - x} \geq \sigma\big(\dist(y, Q)\big)^\mathcal{M}$ for all $ x \in Q $ and $ y \in X $.}
 of order $ \mathcal{M} \geq 1 $, then 
\begin{align} \lab{eq:bh02}
0\leq \dist(y, Q) \leq \left( \dfrac{\emph{Gap}(y, F, X)}{\sigma}\right)^{\frac{1}{\mathcal{M}}}
\quad \mbox{for all} \quad y \in X.
\end{align}
\eprop
\bproo
Let $ y \in \HH $. To prove \eqref{eq:bh}, first note that, by the Cauchy-Schwarz inequality, we have 
$ \inner{H(x)}{y - x} \geq -\norm{H(x)}\norm{y - x}$ for all $ x \in Q $. The result now follows by taking supremum over $ x \in Q $ in both sides of the latter inequality and using the definitions of
$ \mbox{Gap}(y, H, Q) $, $ B_H $ and $ \dist(y, Q) $.
To prove \eqref{eq:bh02}, assume that $ Q $ is $ \sigma $-weakly sharp of order $ \mathcal{M} \geq 1 $. 
Then, for all $ y \in X $, we get 
$ \text{Gap}(y, F, X) \geq \inner{F(x)}{y - x } \geq \sigma\big(\dist(y, Q)\big)^{\mathcal{M}} $ for all 
$ x \in Q $, which clearly gives~\eqref{eq:bh02}.
\eproo

\mgap

\noindent
{\bf Organization of the paper.} 
	In Section \ref{sec:main}, we state our main algorithm, namely IneIREG, the basic assumptions on problem \eqref{eq:prob}\,---\,Assumption \ref{assu:erdos}\,---\,and prove some preliminary results (see Proposition \ref{pro:inv} and Lemma \ref{lem:fd}).
	In Section \ref{sec:hardy}, we present the first results on the iteration-complexity for Algorithm \ref{alg:main} (see Propositions \ref{pro:optm01} and \ref{pro:feas}); specialized results for the case of diminishing and constant regularization parameter sequence
$ \seq{\eta_k} $ are then established in Subsections \ref{subsec:dimish} and \ref{subsec:constant}, respectively.
	Section \ref{sec:strong} is devoted to the study of the strongly monotone case.
	In Section~\ref{sec:ne}, we present some preliminary numerical experiments with three types of problems. Finally, we conclude the paper in Section~\ref{sec:conclusions}.

\section{The main algorithm and some preliminary results} \lab{sec:main}

Consider the bilevel variational inequality problem (VIP) as given in \eqref{eq:prob}, i.e., 
\begin{align} \lab{eq:probmain}
\mbox{Find} \enspace x\in Q \enspace \mbox{such that} \enspace
\Inner{H(x)}{z - x} \geq 0 \quad \text{for all} \enspace z \in Q,
\end{align}
where $ Q $ is the solution set of $ \text{VIP}(F, X) $\,---\,see \eqref{eq:def.Q}\,---\,and the following conditions are assumed to hold:

\bassu \lab{assu:erdos}
\benum
\item \lab{assu:erdos.seg}
$ F \colon \Dom F \subset \HH \to \HH $ and $ H \colon \Dom H \subset \HH \to \HH $ are monotone and Lipschitz continuous maps, with \emph{(}Lipschitz\emph{)} constants $ L_F > 0 $ and $ L_H > 0 $, respectively. 
\item \lab{assu:erdos.ter} 
$ X $ and $ \Omega $ are simple nonempty closed and convex subsets of $ \HH $ such that 
$ X \subset \Omega \subset \Dom F \cap \Dom H$.
\item \lab{assu:erdos.qua}
The solution set $ Q $ of $ \emph{VIP}(F, X) $ is nonempty.
\eenum
\eassu

We propose the following algorithm to numerically solve \eqref{eq:probmain} 
under Assumption \ref{assu:erdos}:

\mgap

\begin{algorithm}[H] \label{alg:main}
\caption{IneIREG: An inertial iteratively regularized EG method for solving \eqref{eq:probmain}}
\SetAlgoLined
\KwInput{Initial guess $ x_0 = x_{-1}\in X $}
\For{$ k = 0, 1, \dots $}{
  Choose an inertial parameter $ \alpha_k > 0 $ and let  
  \begin{align} \lab{eq:def.wk}
   w_k \coloneqq x_k + \alpha_k(x_k - x_{k-1})
  \end{align}
  Choose a stepsize $ \lambda_k > 0 $, a regularization parameter $ \eta_k > 0 $,
  set $ w'_k \coloneqq P_\Omega(w_k) $ and compute 
  \begin{align} \lab{eq:yx}
  \begin{aligned}
  & y _k  = P_X\Big( w_k - \lambda_k\big(F(w'_k) + \eta_kH(w'_k)\big)\Big)\\[2mm]
  & x_{k+1}  =  P_X\Big( w_k - \lambda_k\big(F(y_k) + \eta_kH(y_k)\big)\Big)
  \end{aligned}
  \end{align}
  }
\end{algorithm}

\mgap
\mgap

Next we make some comments regarding Algorithm \ref{alg:main}: 
\begin{enumerate}[label = (\roman*)]
\item As we mentioned before, the operation in \eqref{eq:def.wk}, which generates the (extrapolated) point $ w_k $ from current and past iterate $ x_k $ and $ x_{k-1} $, respectively, introduces \emph{inertial effects} on the iterations generated by Algorithm \ref{alg:main}. Inertial methods, which have their roots in the celebrated heavy-ball and Nesterov's accelerated gradient methods (see, e.g., \cite{Nes18book, Pol87book}), became popular in the context of proximal-type algorithms since the seminal paper~\cite{AlvAtt01}, and were subsequently developed in various research directions (see, e.g., 
\cite{AlvEckGerMel20,AlvMar20,Att21,AttCab20, AttCabChbRia18,BotSedVuo23,CorPey25} and references therein).
\item Note that \eqref{eq:yx}, which represents the main computational burden of Algorithm \ref{alg:main}, is nothing but a step of the celebrated Korpelevich's extragradient (EG) method~\cite{Kor76}
for the monotone operator $ F(\cdot) + \eta_k H(\cdot) $; if we set $ \eta_k \equiv 0 $, then it follows that 
Algorithm \ref{alg:main} reduces to the (inertial) EG method for solving VIP$ (F, X) $. The role of the regularization parameter $ \eta_k >0 $ is to balance, along the iterative process, the two distinct goals while solving \eqref{eq:probmain}, namely optimality (upper-level VIP) and feasibility (lower-level VIP). 
We also mention that the additional projection $ w'_k = P_\Omega(w_k) $ is necessary to recover feasibility of
the extrapolated point $ w_k $ with respect to the domains of 
$ F(\cdot) $ and $ H(\cdot) $ (see Assumption \ref{assu:erdos}\emph{\ref{assu:erdos.ter}}); if 
$ F(\cdot) $ and $ H(\cdot) $ are everywhere defined, then one can simply take $ \Omega = \HH $
and, consequently, $ w'_k \coloneqq w_k $. Finally, we mention that, since the feasible set 
$ X $ is assumed to be simple, it follows that the projection $ P_X(\cdot) $ onto $ X $ can be 
explicitly computed.
\item Conditions on the parameters $ \alpha_k >0 $, $ \lambda_k >0 $ and $ \eta_k >0 $ ensuring theoretical guarantees for Algorithm \ref{alg:main} will be discussed later in this paper; for instance, as it is customary in the analysis of (extra)gradient-type methods, the stepsize $ \lambda_k >0 $ will depend on the inverse of the Lipschitz constant of $ F(\cdot) + \eta_k H(\cdot) $, which happens to be 
$ 1/(L_F + \eta_k L_H) $\,---\, see \eqref{eq:def.Lk} and, e.g., Assumption \ref{assu:var}\emph{\ref{assu:var.seg}} below.
\item Algorithm \ref{alg:main} is a generalization (an inertial version) of 
\cite[Algorithm 3.1]{SamYou25} and many results of this work are inspired by the latter reference.
The first results on the iteration-complexity of Algorithm \ref{alg:main} are Propositions \ref{pro:optm01}
and \ref{pro:feas}, under the Assumption \ref{assu:var}. Theorems \ref{the:main} and \ref{the:ofcons02} (see also \ref{the:ofcons}) present specialized results for
vanishing and constant regularization parameter $ \eta_k > 0 $, respectively. The main results under the assumption that
$ H(\cdot) $ is strongly monotone are presented in Section \ref{sec:strong}.
\end{enumerate}

From now on in this paper, the sequences $\seq{x_k}$,  $ \seq{\alpha_k} $,
 $ \seq{w_k} $, $ \seq{\lambda_k} $, $ \seq{\eta_k} $, $ \seq{w'_k} $ and $ \seq{y_k} $ are generated by
Algorithm \ref{alg:main}. 

\subsection{Preliminary results} \lab{subsec:pre}

In this subsection, we present two preliminary results that will be needed for the theoretical analysis of Algorithm \ref{alg:main}. 
Proposition \ref{pro:inv} is an immediate consequence of the fact (as discussed before) that \eqref{eq:yx} is a step of the EG method for the operator $ F(\cdot) + \eta_k H(\cdot) $; see~\cite{Kor76,MonSva10,Nem05}. 
For the convenience of the reader, we will include a proof of Proposition \ref{pro:inv} in Appendix~B.
On the other hand, Lemma \ref{lem:fd}, which has its roots in~\cite{AlvAtt01}, relates the 
norms of $ w_k $ (extrapolated point) with the norms of current $ x_k $ and past $ x_{k-1} $ iterates.

\bprop \lab{pro:inv}
For all $ x \in X $ and $ k \geq 0 $,
\begin{align*}
%\begin{aligned}
\normq{w_k - x} - \normq{x_{k+1} - x} 
	& \geq 
\Big( 1 - \lambda_k^2\,L_k^2 \Big)\normq{w_k - y_k}\\[2mm]
	&\hspace{2cm} + 2\lambda_k \Inner{F(y_k)}{y_k - x} + 2\lambda_k\eta_k\inner{H(y_k)}{y_k - x},
%\end{aligned}
\end{align*}
where 
\begin{align} \lab{eq:def.Lk}
 L_k \coloneqq L_F + \eta_k L_H  \qquad \mbox{for all} \quad k \geq 0.
\end{align}
\eprop

\blemm \lab{lem:fd}
Let $ x \in X $ and define
\begin{align} \lab{eq:def.pd}
\begin{aligned}
\varphi_k \coloneqq \normq{x_k - x} \quad \text{for all} \quad k\geq -1 \quad \text{and} &\\
	& \hspace{-1cm} \delta_k \coloneqq \alpha_k(1 + \alpha_k)\normq{x_k - x_{k-1}} \quad \text{for all} \quad k\geq 0.
\end{aligned}
\end{align}
Then $ \varphi_{-1} = \varphi_0 $ and
\begin{align} \lab{eq:inervarphi}
\normq{w_k - x} = \varphi_k + \alpha_k(\varphi_k - \varphi_{k-1}) + \delta_k \qquad \mbox{for all} \quad k \geq 0.
\end{align}
\elemm
\bproo
The identity $ \varphi_{-1} = \varphi_0 $ follows from the definition of $ \varphi_k $ 
and the fact that $ x_0 = x_{-1} $ (see the input of Algorithm \ref{alg:main}).
Now, let $ k \geq 0 $ and note that,  from \eqref{eq:def.wk},
\[
 x_k = \dfrac{1}{1+\alpha_k} w_k + \dfrac{\alpha_k}{1+\alpha_k}x_{k-1},
\]
which, in turn, when combined with the well-known identity 
$ \normq{t u + (1-t) v} = t\normq{u} + (1-t)\normq{v} - t(1-t)\normq{u - v}$ yields
\[
\underbrace{\normq{x_k - x}}_{\varphi_k} = \dfrac{1}{1+\alpha_k}\normq{w_k - x} + \dfrac{\alpha_k}{1+\alpha_k}\underbrace{\normq{x_{k-1} - x}}_{\varphi_{k-1}}
- \dfrac{\alpha_k}{(1+\alpha_k)^2}\normq{w_k - x_{k-1}}.
\]
Hence, \eqref{eq:inervarphi} follows simply by multiplying both sides of the above identity by $ 1+ \alpha_k > 0 $ and using the fact (from \eqref{eq:def.wk}) that $ w_k - x_{k-1} = (1+\alpha_k)(x_k - x_{k-1})$, as well as the definition of 
$ \delta_k $ given in \eqref{eq:def.pd}.
\eproo

\section{Iteration-complexity analysis of Algorithm \ref{alg:main}} \lab{sec:hardy}

In this section, we analyze the iteration-complexity of IneIREG. 
In Propositions \ref{pro:optm01} and \ref{pro:feas}, we present preliminary results on \emph{optimality} and \emph{feasibility} for the ergodic mean $ \overline y_k $ (see \eqref{eq:def.ubs} below) with respect to the dual gap functions $ \text{Gap}(\cdot, H, Q) $ and $ \text{Gap}(\cdot, F, X) $, respectively.
In Subsections \ref{subsec:dimish} and \ref{subsec:constant}, we then obtain iteration-complexity 
bounds (optimality and feasibility) for diminishing and constant regularization sequence $ \seq{\eta_k} $, respectively.  
We also refer the reader to Remarks \ref{rem:sensei}, \ref{rem:wets} and \ref{rem:bertsekas} for discussions of our results in the light of the corresponding ones in \cite{SamYou25}.

\mgap

The main results of this section will be established under the following set of assumptions:

\bassu \lab{assu:var}
\benum
\item \lab{assu:var.ter}
The sequence $ \seq{\eta_k} $ is nonincreasing, i.e., $ \eta_k \geq \eta_{k+1} $ for all $ k \geq 0 $.
\item \lab{assu:var.terb}
We have $ \alpha_0 \in [0, 1] $ and the sequence $ \seq{\alpha_k \eta_k^{-1}} $ nonincreasing, i.e., 
$ \alpha_k \eta_k^{-1} \geq \alpha_{k+1} \eta_{k+1}^{-1} $ for all $ k \geq 0 $. 
\item \lab{assu:var.seg}
The sequence of stepsizes $ \seq{\lambda_k} $ is chosen to satisfy: $ \lambda_k \in [\,\underline\lambda, \overline\lambda\,] $ for all $ k \geq 0 $,
where $ 0 < \underline\lambda \leq \overline\lambda \leq 1/L_0 $
and $ L_0 \coloneqq L_F + \eta_0 L_H $.
\item \lab{assu:var.qua} 
The series $ \sum_{k=0}^\infty\,\delta_k \eta_k^{-1} $ is summable, i.e. 
$ s\coloneqq \sum_{k=0}^\infty\, \delta_k \eta_k^{-1} < +\infty $, where $ \delta_k \geq 0 $ is defined as in~\eqref{eq:def.pd}.
\eenum
\eassu

Next we make a few remarks regarding Assumption \ref{assu:var}:

\brema%[\bf On the Assumption \ref{assu:var}] 
\lab{rem:var}
\begin{enumerate}[label = \emph{(\roman*)}]
\item \lab{rem:var.seg}
Note that \emph{Assumptions} \emph{\ref{assu:var}}\ref{assu:var.ter} and 
\emph{\ref{assu:var}}\ref{assu:var.terb} together imply that the sequence of inertial parameters
$ \seq{\alpha_k} $ is also nonincreasing, i.e., $ \alpha_k \geq \alpha_{k+1} $ for all $ k \geq 0 $, 
and $ \alpha_k \in [0, 1] $ for all $ k \geq 0 $.
\item \lab{rem:var.ter}
It results from \emph{Assumptions} \emph{\ref{assu:var}}\ref{assu:var.ter} and 
\emph{\ref{assu:var}}\ref{assu:var.seg}, in particular, that $ 0< \lambda_k \leq 1/L_k $, and so 
$ 1- \lambda_k^2 L_k^2 \geq 0 $ for all $ k \geq 0 $, where $ L_k $ is as in \eqref{eq:def.Lk}.
\item 
Let us now discuss some possible realizations for the sequences $ \seq{\alpha_k } $ and $ \seq{\eta_k } $ 
that ensure the summability condition on \emph{Assumption \ref{assu:var}}\ref{assu:var.qua} holds. Let $ \seq{\eta_k} $ be as in \emph{Assumption \ref{assu:var}}\ref{assu:var.ter} and 
take $\alpha_0 \in [0, 1]$. Then, at each iteration $ k \geq 1 $ of \emph{Algorithm \ref{alg:main}}, one can choose 
$ \alpha_k $ as
\begin{align} \lab{eq:pen}
\alpha_k = 
\begin{cases}
\dfrac{\eta_k}{\eta_{k-1}}\alpha_{k-1},&  \mbox{if} \quad k < m,  \\[3mm]
\eta_k\min\left\{\dfrac{\theta^k}{\normq{x_k - x_{k-1}} + \rho}, \dfrac{\alpha_{k-1}}{\eta_{k-1}}\right\}, &
\mbox{if} \quad k \geq m, 
\end{cases}
\end{align} 
where $ m \in \Z_+ $, $ \rho > 0 $ and $ \theta \in (0, 1) $. 
This choice of $ \seq{\alpha_k} $ clearly ensures that \emph{Assumptions \ref{assu:var}}\ref{assu:var.terb} and 
\emph{\ref{assu:var}}\ref{assu:var.qua} hold.
Note that if $ m = 0 $ in \eqref{eq:pen}, then by a simple argument we have
\begin{align} \lab{eq:bound_s}
	s \leq 2\theta/(1 - \theta), 
\end{align}
where $ s $ is as in \emph{Assumption \ref{assu:var}\emph{\ref{assu:var.qua}}}.
Another possibility is to set $ \eta_k \equiv \eta > 0 $ (constant) and choose $ \alpha_k $, for all $ k \geq 1 $, 
as
\begin{align} \lab{eq:narrow}
\alpha_k = 
\begin{cases}
\xi_0, & \mbox{if} \quad k < m,  \\[3mm]
\xi_{k - m}, &
\mbox{if} \quad k \geq m, 
\end{cases}
\end{align} 
where $ m \in \Z_+$ and $ \seq{\xi_k}_{k\geq 0} $ 
is a nonincreasing summable sequence in $ [0, 1] $. Indeed, first note that, in this case, 
\emph{Assumptions} \emph{\ref{assu:var}}\ref{assu:var.ter} and \emph{\ref{assu:var}}\ref{assu:var.terb} are clearly satisfied and, moreover, \emph{Assumption \ref{assu:var}}\ref{assu:var.qua} reduces to $ \sum_{k=0}^\infty\, \delta_k < + \infty $. 
To conclude the argument, note that from the definition of $ \delta_k $ as in \eqref{eq:def.pd} and the (now) fact that $ \alpha_k \in [0, 1] $, we have, for all $ k \geq 1 $,
$ \delta_k \leq 2\alpha_k \normq{x_k - x_{k-1}} \leq (2D_X^2) \alpha_k $, where in the latter inequality we also used that $ \norm{x_k - x_{k-1}} \leq D_X $ \emph{(}because $ x_k, x_{k-1} \in X $ \emph{)}. Hence, the desired condition follows from the summability of $ \seq{\alpha_k} $ as in \eqref{eq:narrow}.
\item Clearly, by setting $ \alpha_k \equiv 0 $ 
\emph{(}non-inertial version of \emph{Algorithm \ref{alg:main}}\emph{)}, we obtain that \emph{Assumptions}
\emph{\ref{assu:var}}\ref{assu:var.terb} and \emph{\ref{assu:var}}\ref{assu:var.qua} are automatically satisfied.
\end{enumerate}
\erema

\mgap

Next we present the first result on ``optimality'' for Algorithm \ref{alg:main}, in terms of obtaining upper-bounds for the dual gap function $ \text{Gap}(\cdot, H, Q) $ (see Definition \ref{def:dualgap}). 
To this end, we first define, for all $ k \geq 1 $, the \emph{ergodic mean}
\begin{align} \lab{eq:def.ubs}
\overline y_k \coloneqq \dfrac{1}{\Lambda_k}\sum_{j=0}^{k-1}\,\lambda_j y_j \qquad \text{where} \quad
\Lambda_k \coloneqq \sum_{j=0}^{k-1}\lambda_j.
\end{align}
Here, $ \seq{\lambda_k} $ and $ \seq{y_k} $ are generated by Algorithm \ref{alg:main}. Directly from the definition of $ \Lambda_k $ and Assumption \ref{assu:var}\emph{\ref{assu:var.seg}}, we derive that
\begin{align} \lab{eq:play}
	\Lambda_k \geq \underline\lambda\, k \qquad \mbox{for all} \quad k \geq 0.
\end{align}

\mgap

For the rest of this section, recall that $ D_X $ denote the diameter of $ X $.
Also, let $ \seq{\overline y_k} $ and $ \seq{\delta_k} $ be as in \eqref{eq:def.ubs} and \eqref{eq:def.pd}, respectively.

\bprop[\bf Optimality] \lab{pro:optm01}
For all $ k \geq 1 $,
\begin{align} \lab{eq:firstopt}
-B_H \dist(\overline y_k, Q) \leq 
\emph{Gap}(\overline y_k, H, Q) 
	\leq \dfrac{1}{k\,\eta_k}\left(\dfrac{D_X^2}{2\underline \lambda}\right) + 
	\dfrac{\sum_{j=0}^{k-1}\,\delta_j \eta_j^{-1}}{k}\left(\dfrac{1}{2\underline \lambda}\right),
\end{align}
where $ B_H > 0 $ is as in \eqref{eq:def.ch}.
\eprop
\bproo
Let $ x $ be a solution of $ \mbox{VIP}(F, X) $, i.e., $ x $ is an element in  $ Q $ as in \eqref{eq:def.Q}, let 
$ k \geq 1 $ and let $ 0 \leq j \leq k - 1 $. 
Since $ y_j \in X $, it follows from the definition of $ Q $ and the monotonicity of $ F(\cdot) $ that 
$ \inner{F(y_j)}{y_j - x} \geq \inner{F(x)}{y_j - x} \geq 0 $. 
Hence, by using Proposition \ref{pro:inv} and Remark \ref{rem:var}\emph{\ref{rem:var.ter}}, we obtain
\begin{align*}
2\Inner{H(x)}{\lambda_j (y_j - x)} \leq \dfrac{1}{\eta_j}\Big( \normq{w_j - x } - \normq{x_{j+1} - x} \Big),
\end{align*}
where we also used the monotonicity of $ H(\cdot) $ to conclude that 
$ \inner{H(y_j)}{y_j - x} \geq \inner{H(x)}{y_j - x} $.  
 Hence, from Lemma \ref{lem:fd}, the definition of $ \varphi_j $ as in \eqref{eq:def.pd} and some simple algebraic manipulations, we derive that
\begin{align} \lab{eq:print}
\nonumber
2\Inner{H(x)}{\lambda_j (y_j - x)} &\leq \dfrac{1}{\eta_j}\Big( \varphi_j - \varphi_{j+1} +
\alpha_j (\varphi_j - \varphi_{j-1}) + \delta_j \Big)\\[2mm]
%\nonumber
%	& = \dfrac{1}{\eta_j}\left(\varphi_j - \varphi_{j+1}\right) 
%	+ \dfrac{\alpha_j}{\eta_j}\left(\varphi_j - \varphi_{j-1}\right)
%	+ \dfrac{\delta_j}{\eta_j}\\[2mm]
	\nonumber
	& = \dfrac{\varphi_j}{\eta_j} - \dfrac{\varphi_{j+1}}{\eta_{j+1}} + \left(\dfrac{1}{\eta_{j+1}} - \dfrac{1}{\eta_j}\right)\varphi_{j+1} + \dfrac{\alpha_j \varphi_j}{\eta_j} - \dfrac{\alpha_{j-1} \varphi_{j-1}}{\eta_{j-1}}\\[2mm]
	&	\hspace{4cm}+ \left(\dfrac{\alpha_{j-1}}{\eta_{j-1}} - \dfrac{\alpha_j}{\eta_j}\right)\varphi_{j-1}
	+ \dfrac{\delta_j}{\eta_j},
\end{align}
where, for $ j =0 $, we set $ \alpha_{-1} \coloneqq \alpha_0 $ and $ \eta_{-1}\coloneqq \eta_0 $.
Since $ x_{j+1} $, $ x_{j-1} $ and $ x $ belong to  $ X $, we have $ \varphi_{j+1} \leq D_X^2 $
and $ \varphi_{j-1} \leq D_X^2 $. 
Moreover, in view of Assumptions \ref{assu:var}\emph{\ref{assu:var.ter}} and 
\ref{assu:var}\emph{\ref{assu:var.terb}}, we 
also have $ \frac{1}{\eta_{j+1}} - \frac{1}{\eta_j} \geq 0 $ and 
$ \frac{\alpha_{j-1}}{\eta_{j-1}} - \frac{\alpha_j}{\eta_j} \geq 0 $. 
Consequently, from  \eqref{eq:print},
\begin{align} \lab{eq:print02}
\nonumber
2\Inner{H(x)}{\lambda_j (y_j - x)} &\leq \dfrac{\varphi_j}{\eta_j} - \dfrac{\varphi_{j+1}}{\eta_{j+1}}
+ \left(\dfrac{1}{\eta_{j+1}} - \dfrac{1}{\eta_j} + \dfrac{\alpha_{j-1}}{\eta_{j-1}} - \dfrac{\alpha_j}{\eta_j}\right)D_X^2\\[2mm]
& \hspace{3cm }+ \dfrac{\alpha_j \varphi_j}{\eta_j} - \dfrac{\alpha_{j-1} \varphi_{j-1}}{\eta_{j-1}} 
+ \dfrac{\delta_j}{\eta_j}.
\end{align}
	Summing up the inequalities for $ j = 0, \dots, k-1 $ in \eqref{eq:print02}, we then obtain
\begin{align} \lab{eq:zack}
\nonumber
2\Inner{H(x)}{\sum_{j=0}^{k-1}\,\lambda_j(y_j - x)} 
	&\leq \dfrac{\varphi_0}{\eta_0} - \dfrac{\varphi_k}{\eta_k}
+ \left(\dfrac{1}{\eta_k} - \dfrac{1}{\eta_0} + \dfrac{\alpha_0}{\eta_0} - \dfrac{\alpha_{k-1}}{\eta_{k-1}}\right)D_X^2\\[2mm]
\nonumber
	& \hspace{2.5cm} + \left(\dfrac{\alpha_{k-1}}{\eta_{k-1}}\right)\underbrace{\varphi_{k-1}}_{\leq D_X^2} - \left(\dfrac{\alpha_0}{\eta_0}\right)\varphi_0 + \sum_{j=0}^{k-1}\,\dfrac{\delta_j}{\eta_j}\\[2mm]
	\nonumber
	&\leq \underbrace{(1 - \alpha_0)}_{\geq 0}\dfrac{\varphi_0}{\eta_0} 
	+ \left(\dfrac{1}{\eta_k} - \dfrac{1-\alpha_0}{\eta_0}\right)D_X^2 + \sum_{j=0}^{k-1}\,\dfrac{\delta_j}{\eta_j}\\[2mm]
	\nonumber
	&\leq ( 1 - \alpha_0 )\dfrac{D_X^2}{\eta_0} 
	+ \left(\dfrac{1}{\eta_k} - \dfrac{1-\alpha_0}{\eta_0}\right)D_X^2 + \sum_{j=0}^{k-1}\,\dfrac{\delta_j}{\eta_j}\\[2mm]
	& = \dfrac{D_X^2}{\eta_k} + \sum_{j=0}^{k-1}\,\dfrac{\delta_j}{\eta_j},
	%&\leq \left( \dfrac{1}{\eta_k} + \dfrac{\alpha_0}{\eta_0}\right)D_X^2 + s,
\end{align}
where we also used that $ \varphi_0\leq D_X^2 $, $ \varphi_{k-1} \leq D_X^2 $ and $ \alpha_0\in [0, 1] $ (see Assumption \ref{assu:var}\emph{\ref{assu:var.terb}}). Dividing both sides of the inequality \eqref{eq:zack} by $ 2\Lambda_k > 0 $ and using the definitions of $ \overline y_k $ as in \eqref{eq:def.ubs}, we obtain
\begin{align} \lab{eq:zack02}
	\Inner{H(x)}{\overline y_k - x} \leq 
\dfrac{D_X^2}{2\Lambda_k \eta_k} + \dfrac{\sum_{j=0}^{k-1}\, \delta_j\eta_j^{-1}}{2\Lambda_k}.
\end{align}
Now note that the second inequality in \eqref{eq:firstopt} follows directly by
taking supremum in \eqref{eq:zack02} over $ x \in Q $ and using the definition of 
$ \mbox{Gap}(\overline y_k, H, Q) $ (see Definition \ref{def:dualgap}) and the inequality
\eqref{eq:play}. To finish the proof of the proposition, note that the first inequality follows from Proposition~\ref{pro:optm02}.
\eproo

\mgap

\brema
Before passing to the result on feasibility of \emph{IneIREG}, note that there is no guarantee that 
$ \emph{Gap}(\overline y_k, H, Q) \geq 0 $ in \eqref{eq:firstopt}. Indeed, as it was already observed in
\emph{\cite{SamYou25}} for the non-inertial version of \emph{IneIREG}, this is essentially 
because $ \overline y_k $ may be infeasible with respect to the (feasible) set $ Q $ of 
\emph{VIP$ (H, Q) $}; see also \emph{Definition \ref{def:dualgap}}. Following the same 
path of \emph{\cite{SamYou25}}, we will obtain a (negative) lower bound for $ \emph{Gap}(\overline y_k, H, Q) $ under the assumption that the feasible set $ Q $ is 
$ \sigma $-wealky sharp of order $ \mathcal{M} \geq 1 $
\emph{(}see \emph{Proposition \ref{pro:optm02}}, 
and \emph{Theorems \ref{the:main} and \ref{the:ofcons}} below\emph{)}.
\erema

\mgap

\bprop[\bf Feasibility] \lab{pro:feas} 
For all $ k \geq 1 $,
\begin{align} 
0\leq \emph{Gap}(\overline y_k, F, X) 
\lab{eq:4color02}
& \leq \dfrac{1}{k}\left(\dfrac{D_X^2}{2\underline \lambda}\right) 
+ \dfrac{\sum_{j=0}^{k-1}\delta_j}{k}\left(\dfrac{1}{2\underline \lambda}\right)
+ \dfrac{\sum_{j=0}^{k-1}\eta_j}{k} \left(\dfrac{\overline \lambda C_H D_X}{\underline \lambda}\right),
\end{align}
where $ C_H $ is as in \eqref{eq:def.ch}.
\eprop
\bproo
Let $ x \in X $, $ k \geq 1 $ and $ 0\leq j \leq k-1 $. 
Invoking Proposition \ref{pro:inv} and Remark \ref{rem:var}\emph{\ref{rem:var.ter}}, we have
\begin{align} \lab{eq:phone}
\nonumber
2\Inner{F(x)}{\lambda_j(y_j - x)} 
%&\leq 
%\normq{w_j - x} -\normq{x_{j+1} - x} 
%- 2\eta_j\lambda_j\Inner{H(x)}{y_j - x}\\[2mm]
%\nonumber
& \leq \normq{w_j - x} - \underbrace{\normq{x_{j+1} - x}}_{=\varphi_{j+1}} + 2\eta_j\lambda_j\underbrace{\norm{H(x)}}_{\leq C_H}
	\underbrace{\norm{y_j - x}}_{\leq D_X}\\[2mm]
	& \leq \normq{w_j - x} - \varphi_{j+1} + 2\eta_j\lambda_j C_H D_X,
\end{align}
where we also used that 
$ \inner{F(y_j)}{y_j - x} \geq \inner{F(x)}{y_j - x} $ and $ \inner{H(y_j)}{y_j - x} \geq \inner{H(x)}{y_j - x} $
(due to monotonicity of $ F(\cdot) $ and $ H(\cdot) $).
From Lemma \ref{lem:fd}, Remark \ref{rem:var}\emph{\ref{rem:var.seg}} and the fact that 
$ \varphi_{j-1} \leq D_X^2 $, we derive
\begin{align} \lab{eq:phone02}
\nonumber
\normq{w_j - x} - \varphi_{j+1} & = 
%\varphi_j - \varphi_{j+1} + \alpha_j(\varphi_j - \varphi_{j-1}) + \delta_j \\[2mm]
%\nonumber
	%& = 
	\varphi_j - \varphi_{j+1} + \alpha_j \varphi_j - \alpha_{j-1} \varphi_{j-1} 
	+ \underbrace{(\alpha_{j-1} - \alpha_j)}_{\geq 0}\underbrace{\varphi_{j-1}}_{\leq D_X^2} + \delta_j\\[2mm]
	&\leq \varphi_j - \varphi_{j+1} + \alpha_j \varphi_j - \alpha_{j-1} \varphi_{j-1} 
	+ (\alpha_{j-1} - \alpha_j) D_X^2 + \delta_j.
\end{align}
Direct substitution of the upper bound in \eqref{eq:phone02} into \eqref{eq:phone} yields
\begin{align} \lab{eq:phone03} 
\nonumber
2\Inner{F(x)}{\lambda_j(y_j - x)} &\leq \varphi_j - \varphi_{j+1} + \alpha_j \varphi_j - \alpha_{j-1} \varphi_{j-1} 
	+ (\alpha_{j-1} - \alpha_j) D_X^2 + \delta_j \\[2mm]
	& \hspace{2cm} + 2\eta_j\lambda_j C_H D_X.
\end{align}
	Summing up the inequalities for $ j =0, \dots, k-1 $ in \eqref{eq:phone03}, we then obtain
\begin{align} \lab{eq:bourgain}
\nonumber
2\Inner{F(x)}{\sum_{j=0}^{k-1}\lambda_j(y_j - x)} 
	&\leq \varphi_0 - \varphi_k 
+ \alpha_{k-1} \underbrace{\varphi_{k-1}}_{\leq D_X^2} - \alpha_0\varphi_0 + (\alpha_0 - \alpha_{k-1})D_X^2 + \sum_{j=0}^{k-1}\delta_j\\[2mm]
\nonumber
	& \hspace{2cm} + 2\sum_{j=0}^{k-1}\eta_j \left(\overline\lambda C_H D_X\right)  \\[2mm]
	\nonumber
	& \leq \underbrace{(1 - \alpha_0)}_{\geq 0}\underbrace{\varphi_0}_{\leq D_X^2} + \alpha_0 D_X^2 
	+ \sum_{j=0}^{k-1}\delta_j
	+ 2\sum_{j=0}^{k-1}\eta_j \left(\overline\lambda C_H D_X\right)  \\[2mm]
	& \leq D_X^2 + \sum_{j=0}^{k-1}\delta_j 
	+ 2\sum_{j=0}^{k-1}\eta_j \left(\overline\lambda C_H D_X\right),
\end{align}
where, for $ j = 0$, we set $ \alpha_{-1} \coloneqq \alpha_0 $, and
we also used that $ \varphi_{-1} = \varphi_0$ (see Lemma \ref{lem:fd}), 
Assumptions \ref{assu:var}\emph{\ref{assu:var.terb}} and \ref{assu:var}\emph{\ref{assu:var.seg}}
as well as the bounds $ \varphi_{k-1}, \varphi_0 \leq D_X^2 $.
The second inequality in \eqref{eq:4color02} follows directly by dividing both sides of \eqref{eq:bourgain}
by $ 2\Lambda_k >0 $, using the inequality \eqref{eq:play}, the definition of $ \overline y_k $, taking supremum in the resulting inequality over $ x \in X $ and using the definition of $ \text{Gap}(\overline y_k, F, X) $. 
On the other hand, note that the first inequality in 
\eqref{eq:4color02} is a direct consequence of the definition of $ \mbox{Gap}(\overline y_k, F, X) $ and the fact that
$ \overline y_k \in X $ (because it is a convex combination of points in the convex set $X$).
\eproo

\mgap

\brema \lab{rem:sensei}
By setting $ \alpha_k \equiv 0 $, 
we observe that \emph{Propositions \ref{pro:optm01} and \ref{pro:feas}} generalize 
\emph{\cite[Theorem~3.5,(i)–(ii)]{SamYou25}}. Indeed, if $ \alpha_k \equiv 0 $, then it follows that 
$ \delta_k \equiv 0 $, so that \eqref{eq:firstopt} and \eqref{eq:4color02} essentially recover the corresponding results from \emph{\cite{SamYou25}}.
\erema

\subsection{Iteration-complexity of Algorithm \ref{alg:main} with diminishing regularization} \lab{subsec:dimish}

In the next theorem, we summarize the results on optimality and feasibility for the iteration-complexity 
of IneIREG assuming that the following holds:

\bassu \lab{assu:dimish}
The sequence of regularization parameters $ \seq{\eta_k} $ is given by
\begin{align}\lab{eq:def.etak}
 \eta_k \coloneqq \dfrac{\eta_0}{(k+1)^b} \qquad \mbox{for all} \quad k \geq 0,
\end{align}
where $ 0 < \eta_0 \leq 1 $ and $ 0 < b< 1$. 
\eassu 

Clearly, $ \seq{\eta_k} $ as in \eqref{eq:def.etak} satisfies 
Assumption \ref{assu:var}\emph{\ref{assu:var.ter}}.

\btheo[\bf Optimality and Feasibility] \lab{the:main}
Suppose \emph{Assumptions} \emph{\ref{assu:dimish}}, \emph{\ref{assu:var}\emph{\ref{assu:var.terb}}},
\emph{\ref{assu:var}\emph{\ref{assu:var.seg}}} and \emph{\ref{assu:var}\emph{\ref{assu:var.qua}}} hold.
Let $ \seq{\overline y_k} $ be as in 
\eqref{eq:def.ubs}
and let the constants $ C_H $ and $ B_H $ be as in \eqref{eq:def.ch}. 
Then the following statements hold:
\benum
\item \lab{theo:main.seg} 
For all $ k \geq 1 $,
\begin{align} \lab{eq:sadnight}
-B_H \dist(\overline y_k, Q) 
	\leq \emph{Gap}(\overline y_k, H, Q)
 	\leq 
 		\dfrac{1}{k^{1-b}}\left(\dfrac{D_X^2 }{2^{1-b}\underline \lambda\, \eta_0}\right) + 
 		\dfrac{1}{k}\left(\dfrac{s}{2\underline \lambda}\right).
\end{align}
\item \lab{theo:main.ter}
For all $ k \geq 1 $,
\begin{align} \lab{eq:sadnight02} 
0
	\leq \emph{Gap}(\overline y_k, F, X) 
 	\leq \dfrac{1}{k}\left(\dfrac{D_X^2}{2\underline\lambda}\right)
 	+ \dfrac{1}{k}\left(\dfrac{s}{2\underline \lambda}\right)
 		+ \dfrac{1}{k^b}\left(\dfrac{ \eta_0 \overline \lambda C_H D_X}{(1 - b)\underline \lambda}\right).
\end{align}
\item \lab{theo:main.qua}
If $ Q $ is $ \sigma $-weakly sharp of order $ \mathcal{M} \geq 1 $, then, for all $ k \geq 1 $, 
\begin{align} \lab{eq:sadnight03}
 \emph{Gap}(\overline y_k, H, Q) 
 	\geq -\dfrac{B_H}{\sigma^{\,1/\mathcal{M}}} 
 \Bigg( \dfrac{1}{k}\left(\dfrac{D_X^2}{2\underline\lambda}\right)
 	+ \dfrac{1}{k}\left(\dfrac{s}{2\underline \lambda}\right) 
 		+ \dfrac{1}{k^b}\left(\dfrac{ \eta_0 \overline \lambda C_H D_X}{(1 - b)\underline \lambda}\right)
 \Bigg)^
 {\frac{1}{\mathcal{M}}}.
\end{align}
\eenum
\etheo
\bproo
Since $ \seq{\eta_k} $ as in \eqref{eq:def.etak} clearly satisfies Assumption \ref{assu:var}\emph{\ref{assu:var.ter}}, it follows that Assumption \ref{assu:var} hold true. 
This places us in the setting of Propositions \ref{pro:optm01} and \ref{pro:feas}.

\emph{\ref{theo:main.seg}} Note that \eqref{eq:sadnight} follows directly from Proposition \ref{pro:optm01}, the definition of $ s $
as in Assumption \ref{assu:var}\emph{\ref{assu:var.qua}} and the simple inequality 
$ k\eta_k = \eta_0 k(k+1)^{-b} \geq \eta_0 2^{-b} k^{1-b} $  (see \eqref{eq:def.etak}).
%On the other hand, the first inequality in \eqref{eq:sadnight} is a direct consequence of 
%Proposition \ref{pro:optm02}.

\emph{\ref{theo:main.ter}}
We first observe that, in view of \eqref{eq:def.etak},
\begin{align} \lab{eq:integral}
\sum_{j=0}^{k-1}\eta_j = \eta_0\sum_{j=1}^k\,\dfrac{1}{j^b} 
\leq \eta_0 \left(1 + \int_{1}^k\,\dfrac{1}{\tau^b}d\tau\right)
\leq \dfrac{\eta_0}{1-b}k^{1-b}.
\end{align}
Then \eqref{eq:sadnight02} follows from Proposition \ref{pro:feas} combined with \eqref{eq:integral}
and the fact, due to \eqref{eq:def.etak}, that $ \sum_{j=0}^\infty\,\delta_j \leq s $.

\emph{\ref{theo:main.qua}} This is a direct consequence of 
%\eqref{eq:sadnight}, 
\eqref{eq:sadnight02}
and Proposition \ref{pro:optm02}.
\eproo

\mgap

\brema \lab{rem:wets}
Analogously to \emph{Remark \ref{rem:sensei}}, here we also observe that 
\emph{Theorem \ref{the:main}} generalizes \emph{\cite[Theorem 3.6 (Case 1)]{SamYou25}}.
Moreover, the choice of $ \seq{\eta_k} $ as in \eqref{eq:def.etak} also gives
that $ \eta_k/\eta_{k-1} \to 1 $ \emph{(}see \eqref{eq:pen}\emph{)}.
\erema

\subsection{Iteration-complexity of Algorithm \ref{alg:main} for constant regularization} \lab{subsec:constant}

In this subsection, we assume the following:

\bassu \lab{assu:etacons}
The sequence of regularization parameters $ \seq{\eta_k} $ is constant:
\begin{align} \lab{eq:play02}
 \eta_k \equiv \eta > 0 \qquad \mbox{for all} \quad k \geq 0.
\end{align}
\eassu

\mgap

\brema[\bf On the Assumptions \ref{assu:etacons} and \ref{assu:var}] \lab{rem:block}
We notice that \eqref{eq:play02} implies that \emph{Assumption \ref{assu:var}\emph{\ref{assu:var.ter}}} is obviously satisfied and \emph{Assumption \ref{assu:var}\emph{\ref{assu:var.terb}}}
reduces to $ \seq{\alpha_k} $ being nonincreasing, i.e., $ \alpha_k \geq \alpha_{k+1} $ for all
$ k \geq 0 $. 
Moreover,  \emph{Assumption \ref{assu:var}\emph{\ref{assu:var.qua}}} now reduces to the serie
$ \sum_{k=0}^\infty \delta_k $ being summable, i.e.,
\begin{align} \lab{eq:def.heta}
 \what s \coloneqq \sum_{k=0}^\infty \delta_k < + \infty.
\end{align}
\erema

\btheo[\bf Optimality and Feasibility] \lab{the:ofcons}
Suppose \emph{Assumptions} \emph{\ref{assu:etacons}}, \emph{\ref{assu:var}\emph{\ref{assu:var.terb}}},
\emph{\ref{assu:var}\emph{\ref{assu:var.seg}}} and \emph{\ref{assu:var}\emph{\ref{assu:var.qua}}} hold\,---\,see \emph{Remark \ref{rem:block}}.
Let $ \seq{\overline y_k} $ be as in \eqref{eq:def.ubs}, let $ \what s $ be as in \eqref{eq:def.heta} and let 
$ B_H $ and $ C_H $ be as in \eqref{eq:def.ch}. 
Then the following statements hold: 
\benum
\item \lab{the:ofcons.seg} 
For all $ k \geq 1 $,
\begin{align} \lab{eq:sadday}
-B_H \dist(\overline y_k, Q) 
	\leq \emph{Gap}(\overline y_k, H, Q)
 	\leq 
 	\dfrac{1}{k}\left( \dfrac{D_X^2 + \what s}{2\underline\lambda \eta} \right).
\end{align}
\item \lab{the:ofcons.ter}
For all $ k \geq 1 $,
\begin{align} \lab{eq:sadday02}
0 
	\leq \emph{Gap}(\overline y_k, F, X)
 	\leq \dfrac{1}{k}\left( \dfrac{D_X^2 + \what s}{2\underline\lambda} \right)
		+ \eta \left( \dfrac{\overline\lambda C_H D_X}{\underline\lambda}\right).
\end{align}
\item \lab{the:ofcons.qua}
If $ Q $ is $ \sigma $-weakly sharp of order $ \mathcal{M} \geq 1 $, then, for all $ k \geq 1 $, 
\begin{align} \lab{eq:sadday03}
 \emph{Gap}(\overline y_k, H, Q) 
 	\geq -\dfrac{B_H}{\sigma^{\,1/\mathcal{M}}} 
 		\Bigg( \dfrac{1}{k}\left( \dfrac{D_X^2 + \what s}{2\underline\lambda} \right)
		+ \eta \left( \dfrac{\overline\lambda C_H D_X}{\underline\lambda}\right) \Bigg)^{\frac{1}{\mathcal{M}}}.
\end{align}
\eenum
\etheo
\bproo
Analogously to Theorem \ref{the:main}'s proof, we are again in the setting of Propositions \ref{pro:optm01}
and \ref{pro:feas}\,---\,see Remark \ref{rem:block}.

\emph{\ref{the:ofcons.seg}} This follows from Proposition \ref{pro:optm01},
Assumption \ref{assu:etacons} and the definition of $ \what s $ as in \eqref{eq:def.heta}.

\emph{\ref{the:ofcons.ter}} This is a direct consequence of Proposition \ref{pro:feas},
Assumption \ref{assu:etacons} and the definition of $ \what s $.

\emph{\ref{the:ofcons.qua}} This follows from \eqref{eq:sadday02}
and Proposition \ref{pro:optm02}.
\eproo

\mgap

\btheo[\bf Optimality and Feasibility] \lab{the:ofcons02}
Suppose all assumptions of \emph{Theorem \ref{the:ofcons}} hold and, additionally, that
\begin{align} \lab{eq:etaconsb}
 \eta \coloneqq \dfrac{\veps}{2\mathcal{D}_0},
\end{align}
where $ \veps >0 $ is a given tolerance and $ \mathcal{D}_0 > 0 $ is such that
\begin{align} \lab{eq:ubdz}
\mathcal{D}_0 \geq \dfrac{\overline \lambda C_H D_X}{\underline \lambda}.
\end{align}
Then, for all
\begin{align} \lab{eq:lagrange}
k \geq 
\max\left\{ 
\left\lceil \dfrac{1}{\veps^2}\left( \dfrac{\mathcal{D}_0(D_X^2 + \what s)}{\underline\lambda}\right) \right \rceil,
\left\lceil \dfrac{1}{\veps}\left( \dfrac{D_X^2 + \what s}{\underline\lambda}\right)\right\rceil
\right\},
\end{align}
we have
\begin{align} \lab{eq:bertsekas}
-B_H \dist(\overline y_k, Q) \leq  
\emph{Gap}(\overline y_k, H, Q) \leq \veps \quad \mbox{and} \quad 
0\leq \emph{Gap}(\overline y_k, F, X) \leq \veps.
\end{align}
Moreover, if $ Q $ is $ \sigma $-weakly sharp of order $ \mathcal{M} \geq 1 $, then 
\begin{align} \lab{eq:bertsekas02}
\emph{Gap}(\overline y_k, H, Q) \geq -\left(\dfrac{B_H}{\sigma^{1/\mathcal{M}}}\right)
 \veps^{\frac{1}{\mathcal{M}}}.
\end{align}
\etheo
\bproo
The proof follows readily from Theorem \ref{the:ofcons}, \eqref{eq:etaconsb} and \eqref{eq:ubdz}.
\eproo

\mgap

\brema \lab{rem:bertsekas}
In \emph{\cite{SamYou25}},  the case of constant regularization parameter was analyzed in \emph{Theorem 3.6 (Case 2)}. 
We mention that while our bound on optimality and feasibility given in \eqref{eq:lagrange} and \eqref{eq:bertsekas} is independent on assumptions of
$ Q $ being $ \sigma $-weakly sharp of any order $ \mathcal{M} \geq 1 $,  the corresponding results in 
\emph{Theorem 3.6 (Case 2)} are entirely supported on the assumption that $ Q $ is 
$ \sigma $-weakly sharp of order $ \mathcal{M} = 1 $. Moreover, we also mention that the complexity analysis in the latter theorem
(Case 2) also depends on the assumption\footnote{In the notation of \cite{SamYou25}, $ \alpha $ is used instead of $ \sigma $.} $\eta \leq \frac{\sigma}{2\norm{H(x^*)}}$, which is similar to our assumption 
\eqref{eq:ubdz}.
\erema

\section{Iteration-complexity analysis of Algorithm \ref{alg:main} under strong monotonicity} \lab{sec:strong}

In this section, additionally to Assumption \ref{assu:erdos}, we assume 
the following about Problem \eqref{eq:probmain}:

\bassu \lab{assu:hmu}
The operator $ H(\cdot) $ is \emph{$ \mu $-strongly monotone}, i.e., there exists $ \mu > 0 $ such that
\begin{align} \lab{eq:hmu}
	\inner{H(x) - H(y)}{x - y} \geq \mu\normq{x - y} \qquad \mbox{for all} \quad x, y \in \Dom H.
\end{align}
\eassu

In Propositions \ref{pro:optstr} and \ref{pro:feasstr}, we will present the first results on convergence rates for optimality and feasibility, respectively. These results will be then refined in Subsection \ref{subsec:str} for the case where the regularization parameter sequence $ \seq{\eta_k} $ is assumed to be constant.

\mgap

In the next proposition, we use Assumption \ref{assu:hmu} to improve the result of Proposition \ref{pro:inv}:

\bprop \lab{pro:invstr}
For all $ x \in X $ and $ k \geq 0 $,
\begin{align} \lab{eq:happynow}
\begin{aligned}
\normq{w_k - x} - \normq{x_{k+1} - x} & \geq 
\Big( 1 - \lambda_k^2 L_k^2 \Big)\normq{w_k - y_k} + 2\lambda_k \Inner{F(x)}{y_k - x} \\[2mm]
&\hspace{2cm} + 2\lambda_k\eta_k\inner{H(x)}{y_k - x} 
+ 2\lambda_k\eta_k\mu\normq{y_k - x},
\end{aligned}
\end{align}
where $ L_k > 0 $ is as in \eqref{eq:def.Lk}.
\eprop
\bproo
Let $ x \in X $ and $ k \geq 0 $. From the monotonicity of $ F(\cdot) $ and \eqref{eq:hmu}, we have $ \inner{F(y_k)}{y_k - x} \geq \inner{F(x)}{y_k - x} $ and
$ \inner{H(y_k)}{y_k - x} \geq \inner{H(x)}{y_k - x} + \mu\normq{y _k - x} $, respectively, which in turn gives
\eqref{eq:happynow} as a direct consequence of Proposition \ref{pro:inv}.
\eproo

The iteration-complexity analysis of Algorithm \ref{alg:main} will also rely, in particular, on the following set of assumptions on $ \seq{\lambda_k} $, $ \seq{\alpha_k} $ and $ \seq{\eta_k} $:

\bassu \lab{assu:lorenz}
\benum
\item \lab{assu:lorenz.ter}
The sequence of stepsizes $ \seq{\lambda_k} $ is chosen to satisfy:
$0 < \lambda_k < 1/L_k $, for all $ k \geq 0 $, 
where $ L_k > 0 $ is as in \eqref{eq:def.Lk}.
\item \lab{assu:lorenz.qua}
We have $ \alpha_0 \in [0, 1] $ and $ \alpha_{k+1}\leq ( 1- \beta_k)\alpha_k $ for all $ k \geq 0 $, where 
\begin{align} \lab{eq:def.betak}
\beta_k \coloneqq \left( \dfrac{1}{1-\lambda_k^2 L_k^2} + \dfrac{1}{2\lambda_k \eta_k \mu} \right)^{-1} \in (0, 1).
\end{align}
\eenum
\eassu

\mgap

\brema \lab{rem:osher}
In \emph{\cite{SamYou25}}, the algorithm \emph{IR-EG}$_{\verb|m,m|}$ is 
accordingly modified and tailored for the strongly monotone case 
\emph{(}see \emph{IR-EG}$_{\verb|s,m|}$\emph{)}. In particular, the stepsize $ \gamma > 0 $ 
which \emph{(}in the notation of \emph{\cite{SamYou25}}\emph{)} was chosen to satisfy
$ \gamma^2(L_F^2 + \eta_0 L_H^2) \leq 0.5 $ \emph{(}see \emph{\cite[Theorem 3.5]{SamYou25}}\emph{)} 
was then modified to satisfy the more restrictive condition $ \gamma^2L_F^2 + \gamma\eta_k\mu_H + \gamma^2\eta_k^2L_H^2 \leq 0.5 $ \emph{(}see \emph{IR-EG}$_{\verb|s,m|}$\emph{)}.
We, on the other hand, do not modify our main algorithm (\emph{Algorithm \ref{alg:main}}) for the strongly monotone case, in this way keeping the same stepsize assumptions for both the monotone and strongly monotone problems; see \emph{Assumptions \ref{assu:var}\emph{\ref{assu:var.seg}}, \ref{assu:lorenz}\emph{\ref{assu:lorenz.ter}} and \ref{assu:botelho}\emph{\ref{assu:botelho.ter}}(below)}.
This was made possible thanks to an optimization in \emph{Corollary \ref{cor:invstr}}'s proof produced by employing \emph{(}the simple\emph{)} \emph{Lemma \ref{lem:min}} given in Appendix~A.
\erema

\mgap

\bcoro \lab{cor:invstr}\emph{[of Proposition \ref{pro:invstr}]}
Suppose \emph{Assumption \ref{assu:lorenz}\emph{\ref{assu:lorenz.ter}}} holds 
and let $ \seq{\beta_k} $ as in \eqref{eq:def.betak}.
Then, for all $ x \in X $ and $ k\geq 0 $,
\begin{align} \lab{eq:happynow02}
(1 - \beta_k)\normq{w_k - x} - \normq{x_{k+1} - x} \geq 2\lambda_k\inner{F(x)}{y_k - x} +
2\lambda_k \eta_k \inner{H(x)}{y_k - x}.
\end{align}
\ecoro
\bproo
Let $ x \in X $ and $ k \geq 0 $. In view of Lemma \ref{lem:min} (see Appendix~A) with $ a \coloneqq 1 - \lambda_k^2L_k^2 $, 
$ b \coloneqq 2\lambda_k\eta_k\mu $ and $ c \coloneqq \norm{w_k - x} $, and the definition of 
$ \beta_k \in (0, 1) $, we get
\begin{align} \lab{eq:lee}
\left(1 - \lambda_k^2L_k^2\right)\normq{w_k - y_k} + 2\lambda_k\eta_k\mu\normq{y_k - x} \geq 
\beta_k \normq{w_k - x},
\end{align}
where we also used the Cauchy-Schwarz inequality to conclude that  
$ \norm{w_k - y_k} + \norm{y_k - x} \geq \norm{w_k - x} $, which, in turn, 
yields \eqref{eq:happynow02} from \eqref{eq:lee}
and Proposition \ref{pro:invstr}.
\eproo

\mgap

From now on in this section, unless stated otherwise, we suppose Assumption \ref{assu:lorenz} holds.
Define
\begin{align} \lab{eq:def.pk}
p_{-1} \coloneqq 1 \quad \mbox{and} \quad 
p_k \coloneqq \left( \prod_{i=0}^k\,(1 - \beta_i)\right)^{-1} \qquad \mbox{for all}\quad k \geq 0,
\end{align}
where $ \seq{\beta_k} $ is as in \eqref{eq:def.betak}. 
Directly from the definition of $ \seq{p_k }$ and
Assumption \ref{assu:lorenz}\emph{\ref{assu:lorenz.qua}}, we obtain that 
the sequence $ \seq{p_{k-1}\alpha_k} $ is nonincreasing, i.e.,
\begin{align} \lab{eq:pkapk}
 p_{k-1}\alpha_k - p_k\alpha_{k+1}\geq 0 \qquad \mbox{for all}\quad k \geq 0.
\end{align}
Indeed, note that, in this case, we have 
$ p_{k-1}\alpha_k - p_k\alpha_{k+1} \geq  p_{k-1}\alpha_k - \underbrace{p_k(1 - \beta_k)}_{p_{k-1}}\alpha_k = 0 $,
for all $ k \geq 0 $.

\mgap

Similarly to \eqref{eq:def.ubs}, we now define the sequence of \emph{ergodic means} 
$ \seq{\overline y_k} $, for all $ k \geq 1 $,
by
\begin{align} \lab{eq:def.ubs02}
\overline y_k \coloneqq \dfrac{1}{\Lambda_k}\sum_{j=0}^{k-1}\,(\lambda_j \eta_j p_j) y_j\,, \qquad \text{where} \quad
\Lambda_k \coloneqq \sum_{j=0}^{k-1}\lambda_j \eta_j p_j\,,
\end{align}
where $ \seq{p_k} $ is as in \eqref{eq:def.pk}.

\mgap

For the remainder of this paper, recall that $ D_X $ denotes the diameter of $ X $.
Also, let $ \seq{\overline y_k} $ and $ \seq{\Lambda_k} $ be as in \eqref{eq:def.ubs02} and
let $ \seq{p_k} $ and $ \seq{\delta_k} $ be as in \eqref{eq:def.pk} and \eqref{eq:def.pd}, respectively.

\mgap

\bprop[\bf Optimality] \lab{pro:optstr}
For all $ k \geq 1 $,
\begin{align} \lab{eq:doc}
-B_H \dist(\overline y_k, Q) \leq 
\emph{Gap}(\overline y_k, H, Q)
	\leq \dfrac{1}{\Lambda_k}\left( \dfrac{D_X^2 + \sum_{j=0}^{k-1}\,p_{j-1}\delta_j}{2} \right).
\end{align}
\eprop
\bproo
Let $ k \geq 1 $ and $ 0 \leq j \leq k -1 $. 
Also, let $ x \in Q $, i.e., let $ x \in X $ be a solution of $ \text{VIP}(F, X) $. Then, since in this case
$ \inner{F(x)}{y_j - x} \geq 0 $, from Corollary \ref{cor:invstr}, we obtain
\begin{align*}
(1 - \beta_j)\normq{w_j - x} - \normq{x_{j+1} - x} \geq 2\lambda_j \eta_j \inner{H(x)}{y_j - x}.
\end{align*}
Multiplying both sides of the latter equation by $ p_j $, we get
\begin{align*}
\underbrace{p_j(1 - \beta_j)}_{p_{j-1}}\normq{w_j - x} - p_j\normq{x_{j+1} - x} \geq 2\lambda_j \eta_j p_j\inner{H(x)}{y_j - x},
\end{align*}
which combined with Lemma \ref{lem:fd} and the definition of $ \varphi_{j+1} $ yields
\begin{align} \lab{eq:eliz}
 p_{j-1}\bigg(\varphi_j + \alpha_j(\varphi_j - \varphi_{j-1}) + \delta_j \bigg) - p_j\varphi_{j+1} \geq 
 2\lambda_j \eta_j p_j\inner{H(x)}{y_j - x}.
\end{align}
Now note that, after some simple algebraic manipulations,
\begin{align} \lab{eq:buarque}
\nonumber
& p_{j-1}\bigg(\varphi_j + \alpha_j(\varphi_j - \varphi_{j-1}) + \delta_j \bigg) - p_j\varphi_{j+1} \\[2mm]
%\nonumber
%& = p_{j-1}\varphi_j - p_j \varphi_{j+1} + p_{j-1}\alpha_j (\varphi_j -  \varphi_{j-1}) + p_{j-1}\delta_j \\[2mm]
%\nonumber
%& = p_{j-1}\varphi_j - p_j \varphi_{j+1} + p_{j-1}\alpha_j\varphi_j - p_{j-1}\alpha_j \varphi_{j-1} + p_{j-1}\delta_j 
%\\[2mm]
& = p_{j-1}\varphi_j - p_j \varphi_{j+1} + p_j\alpha_{j+1}\varphi_j - p_{j-1}\alpha_j \varphi_{j-1} 
+ \left(p_{j-1}\alpha_j - p_j\alpha_{j+1}\right)\varphi_j + p_{j-1}\delta_j.
\end{align}
From \eqref{eq:pkapk}, \eqref{eq:buarque} and the fact that $ \varphi_k \leq D_X^2 $, we have
\begin{align*}
p_{j-1}\bigg(\varphi_j + \alpha_j(\varphi_j - \varphi_{j-1}) + \delta_j \bigg) - p_j\varphi_{j+1} & \leq
p_{j-1}\varphi_j - p_j \varphi_{j+1} + p_j\alpha_{j+1}\varphi_j - p_{j-1}\alpha_j \varphi_{j-1}\\[2mm] 
& \hspace{1.5cm} + \left(p_{j-1}\alpha_j - p_j\alpha_{j+1}\right) D_X^2 + p_{j-1}\delta_j,
\end{align*}
which, in turn, combined with \eqref{eq:eliz} yields
\begin{align*}
2\lambda_j \eta_j p_j\inner{H(x)}{y_j - x} &\leq 
p_{j-1}\varphi_j - p_j \varphi_{j+1} + p_j\alpha_{j+1}\varphi_j - p_{j-1}\alpha_j \varphi_{j-1}\\[2mm]
& \hspace{1.5cm} + \left(p_{j-1}\alpha_j - p_j\alpha_{j+1}\right) D_X^2 + p_{j-1}\delta_j.
\end{align*}
Summing up the latter inequality for $ j = 0, \dots, k -1 $, we derive
\begin{align} \lab{eq:eliz02}
\nonumber
2\Inner{H(x)}{\sum_{j=0}^{k-1}\,\lambda_j \eta_j p_j(y_j - x)} 
	&\leq \underbrace{p_{-1}}_{=1}\varphi_0 - p_{k-1}\varphi_k + p_{k-1}\alpha_k \underbrace{\varphi_{k-1}}_{\leq D_X^2} - 
\underbrace{p_{-1}}_{=1}\alpha_0\underbrace{\varphi_{-1}}_{=\varphi_0}\\[2mm]
\nonumber
	&\hspace{1cm} + \big(\underbrace{p_{-1}}_{=1}\alpha_0 - p_{k-1}\alpha_k\big) D_X^2 + \sum_{j=0}^{k-1}\,p_{j-1}\delta_j\\[2mm]
\nonumber
	& \leq \underbrace{(1- \alpha_0)}_{\geq 0}\underbrace{\varphi_0}_{\leq D_X^2} + \alpha_0 D_X^2 + \sum_{j=0}^{k-1}\,p_{j-1}\delta_j\\[2mm]
& \leq D_X^2 + \sum_{j=0}^{k-1}\,p_{j-1}\delta_j.
\end{align}
By dividing both sides of \eqref{eq:eliz02} by $ 2\Lambda_k > 0 $ and using the definition of $ \overline y_k $, we then get
\begin{align*}
\Inner{H(x)}{\overline y_k - x} \leq \dfrac{1}{\Lambda_k}
\left( \dfrac{D_X^2 + \sum_{j=0}^{k-1}\,p_{j-1}\delta_j}{2} \right).
\end{align*}
Then the second inequality in \eqref{eq:doc} follows by taking the supremum in the latter inequality over $ x \in Q $ and using the definition of $ \text{Gap}(\overline y_k, H, Q) $. To finish the proof of the proposition, note that the first inequality is a direct consequence of Proposition \ref{pro:optm02}.
\eproo

\mgap

\bprop[\bf Feasibility] \lab{pro:feasstr}
Suppose, additionally to \emph{Assumption \ref{assu:lorenz}}, that 
\emph{Assumption} \emph{\ref{assu:var}\emph{\ref{assu:var.ter}}} holds on $ \seq{\eta_k} $, i.e.,
suppose $ \seq{\eta_k} $ is nonincreasing.
Then, for all $ k \geq 1 $,
\begin{align} \lab{eq:doc02}
0 \leq \emph{Gap}(\overline y_k, F, X)
	\leq 
\dfrac{1}{2\Lambda_k}\left( \eta_0 D_X^2 
+ \sum_{j=0}^{k-1}\,\eta_jp_{j-1}\delta_j 
+ 2\left(\sum_{j=0}^{k-1}\,\lambda_j \eta_j^2 p_j\right) C_H D_X
\right),
\end{align}
where $ C_H $ is as in \eqref{eq:def.ch}.
\eprop
\bproo
Let $ k \geq 1 $, let $ 0 \leq j \leq k-1 $ and let $ x \in X $. From Corollary \ref{cor:invstr} and the Cauchy-Schwarz inequality, we have
\begin{align*}
2\lambda_j \inner{F(x)}{y_j - x} & \leq (1 - \beta_j)\normq{w_j - x} - \normq{x_{j+1} - x} 
- 2\lambda_j \eta_j \inner{H(x)}{y_j - x}\\[2mm]
& \leq (1 - \beta_j)\normq{w_j - x} - \normq{x_{j+1} - x} 
+ 2\lambda_j \eta_j \underbrace{\norm{H(x)}}_{\leq C_H}\underbrace{\norm{y_j - x}}_{\leq D_X}.
\end{align*}
By multiplying both sides of the latter inequality by $ \eta_j p_j \geq 0 $, we find
\begin{align} \lab{eq:power}
2\lambda_j \eta_j p_j \inner{F(x)}{y_j - x} \leq 
\eta_j \underbrace{p_j (1 - \beta_j)}_{p_{j-1}}\normq{w_j - x} - \eta_j p_j \normq{x_{j+1} - x}
+ 2\left(\lambda_j \eta_j^2 p_j\right) C_H D_X.
\end{align}
Now note that, from Lemma \ref{lem:fd} and Assumption \ref{assu:var}\emph{\ref{assu:var.ter}} (i.e., $\eta_j\geq \eta_{j+1}$), we obtain
\begin{align} \lab{eq:again00}
\nonumber
	& \eta_j p_{j-1}\normq{w_j - x} - \eta_j p_j \normq{x_{j+1} - x}\\[2mm]
%\nonumber
%	& = \eta_j p_{j-1}\Big( \varphi_j + \alpha_j(\varphi_j - \varphi_{j-1}) + \delta_j\Big) - \eta_j p_j 
%\varphi_{j+1}\\[2mm]
\nonumber
	& \leq \eta_j p_{j-1}\Big( \varphi_j + \alpha_j(\varphi_j - \varphi_{j-1}) + \delta_j\Big) 
    - \eta_{j+1} p_j \varphi_{j+1} \\[2mm]
	& = \eta_j p_{j-1}\varphi_j - \eta_{j+1} p_j \varphi_{j+1} + \eta_j p_{j-1}\alpha_j(\varphi_j - \varphi_{j-1})
+ \eta_j p_{j-1}\delta_j.
\end{align}
Note also that
\begin{align} \lab{eq:again0}
\nonumber
 \eta_j p_{j-1}\alpha_j(\varphi_j - \varphi_{j-1}) & = 
 \eta_j p_{j-1}\alpha_j\varphi_j - \eta_j p_{j-1}\alpha_j\varphi_{j-1}\\[2mm]
	 & = \eta_{j+1} p_j\alpha_{j+1}\varphi_j - \eta_j p_{j-1}\alpha_j\varphi_{j-1} +
 \Big(\eta_j p_{j-1}\alpha_j - \eta_{j+1} p_j\alpha_{j+1}\Big)\underbrace{\varphi_j}_{\leq D_X^2}.
\end{align}
Now we observe that, since $ \eta_j \geq \eta_{j+1} $, we have
\begin{align} \lab{eq:again}
\nonumber
\eta_j p_{j-1}\alpha_j - \eta_{j+1} p_j\alpha_{j+1} & \geq \eta_{j+1} p_{j-1}\alpha_j - \eta_{j+1} p_j\alpha_{j+1}\\[2mm]
\nonumber
	& = \eta_{j+1} \underbrace{\Big(p_{j-1}\alpha_j - p_j\alpha_{j+1} \Big)}_{\geq 0}\\[2mm]
	& \geq  0,
\end{align}
where in the latter inequality we used condition \eqref{eq:pkapk}.
By combining \eqref{eq:again00}--\eqref{eq:again}, we get
\begin{align*}
&\eta_j p_{j-1}\normq{w_j - x} - \eta_j p_j \normq{x_{j+1} - x} \\[2mm] 
	& \leq \eta_j p_{j-1}\varphi_j - \eta_{j+1} p_j \varphi_{j+1} +
\eta_{j+1} p_j\alpha_{j+1}\varphi_j - \eta_j p_{j-1}\alpha_j\varphi_{j-1} \\[2mm] 
	& \hspace{5cm}+ \Big(\eta_j p_{j-1}\alpha_j - \eta_{j+1} p_j\alpha_{j+1}\Big) D_X^2
+ \eta_j p_{j-1}\delta_j
\end{align*}
and, from \eqref{eq:power}, 
\begin{align} \lab{eq:power02}
\nonumber
&2\lambda_j \eta_j p_j \inner{F(x)}{y_j - x} \\[2mm] 
%\nonumber
%&\leq \eta_j p_{j-1}\normq{w_j - x} - \eta_j p_j \normq{x_{j+1} - x}
%+ 2\left(\lambda_j \eta_j^2 p_j\right) C_H D_X \\[2mm]
\nonumber
	& \leq \eta_j p_{j-1}\varphi_j - \eta_{j+1} p_j \varphi_{j+1} +
\eta_{j+1} p_j\alpha_{j+1}\varphi_j - \eta_j p_{j-1}\alpha_j\varphi_{j-1} \\[2mm] 
	& \hspace{2cm}+ \Big(\eta_j p_{j-1}\alpha_j - \eta_{j+1} p_j\alpha_{j+1}\Big) D_X^2
+ \eta_j p_{j-1}\delta_j + 2\left(\lambda_j \eta_j^2 p_j\right) C_H D_X.
\end{align}
Summing up \eqref{eq:power02} for $ j = 0, \dots, k - 1 $, we find
\begin{align} \lab{eq:paper}
\nonumber
	& 2\Inner{F(x)}{\sum_{j=0}^{k-1}\lambda_j \eta_j p_j(y_j - x)}\\[2mm]
\nonumber
	& \leq \eta_0 \underbrace{p_{-1}}_{=1} \varphi_0 - \eta_k p_{k-1}\varphi_k
+ \eta_k p_{k-1}\alpha_k \underbrace{\varphi_{k-1}}_{\leq D_X^2} 
- \eta_0 \underbrace{p_{-1}}_{=1}\alpha_0 \underbrace{\varphi_{-1}}_{=\varphi_0}\\[2mm]
\nonumber
	& \hspace{2cm}+ \Big(\eta_0 \underbrace{p_{-1}}_{=1}\alpha_0 - \eta_k p_{k-1}\alpha_k\Big) D_X^2
+ \sum_{j=0}^{k-1}\,\eta_j p_{j-1}\delta_j + 2\left(\sum_{j=0}^{k-1}\lambda_j \eta_j^2 p_j\right) C_H D_X\\[2mm]
\nonumber
	& \le \underbrace{(1 - \alpha_0)}_{\geq 0}\eta_0\underbrace{\varphi_0}_{\leq D_X^2} + \eta_0\alpha_0 D_X^2
+ \sum_{j=0}^{k-1}\,\eta_j p_{j-1}\delta_j + 2\left(\sum_{j=0}^{k-1}\lambda_j \eta_j^2 p_j\right) C_H D_X \\[2mm]
	& \leq \eta_0 D_X^2 + \sum_{j=0}^{k-1}\,\eta_j p_{j-1}\delta_j + 2\left(\sum_{j=0}^{k-1}\lambda_j \eta_j^2 p_j\right) C_H D_X.
\end{align}
To finish the proof of the proposition, note that the second inequality in \eqref{eq:doc02} follows directly by dividing both sides of \eqref{eq:paper} by $ 2\Lambda_k > 0 $, using the definition of $ \overline y_k $ and, subsequently, taking supremum of the resulting inequality over all $ x \in X $, and using the definition of 
$ \text{Gap}( \overline y_k, F, X ) $. Recall that the first inequality in \eqref{eq:doc02} is a direct consequence of the fact that $ \overline y_k \in  X $ for all $ k \geq 1 $.
\eproo

\subsection{Iteration-complexity of Algorithm \ref{alg:main} (under strong monotonicity) for constant regularization} 
\lab{subsec:str}

In this subsection, additionally to Assumption \ref{assu:lorenz}, unless stated otherwise, we make the following assumptions on the execution of Algorithm \ref{alg:main} (cf. Assumption \ref{assu:var}):

\bassu \lab{assu:botelho}
\benum
\item \lab{assu:botelho.seg}
The sequence of regularization parameters $ \seq{\eta_k} $ is constant, i.e.,\\
$ \eta_k \equiv \eta > 0 $ for all $ k \geq 0 $.
\item \lab{assu:botelho.ter}
The sequence of stepsizes $ \seq{\lambda_k} $ is chosen to satisfy:
$ \lambda_k \in [\,\underline\lambda, \overline\lambda\,] $ for all $ k \geq 0 $,
where $ 0 < \underline\lambda \leq \overline\lambda < 1/L $
and $ L \coloneqq L_F + \eta L_H $.
\eenum
\eassu

\mgap

The main results on optimality and feasibility are given in Corollaries \ref{cor:optstr02} and 
\ref{cor:feasstr02}, respectively (see also Corollaries \ref{cor:optstr03} and \ref{cor:feasstr03}). 
The main result of this subsection is then presented in Theorem \ref{the:ejinho}. 
We also refer the reader to Remarks \ref{rem:rockafellar}, \ref{rem:shor} and \ref{rem:suasuna} for additional discussions.

\mgap

We begin with the following (technical) lemma:

\blemm \lab{lem:raman}
For all $ k \geq 1 $,
\begin{align} \lab{eq:belchior}
\Lambda_k 
	\geq \underline \lambda \eta \sum_{j=0}^{k-1}p_j 
	\geq \dfrac{\underline\lambda \eta }{(1-\beta)^k},
\end{align}
where 
\begin{align} \lab{eq:def.beta}
\beta \coloneqq \left(\dfrac{1}{1 - \overline\lambda^{\,2} L^2} 
+ \dfrac{1}{2\underline\lambda \eta \mu}\right)^{-1} \in (0, 1).
\end{align}
\elemm
\bproo
Let $ \seq{\beta_k} $ be as in \eqref{eq:def.betak}; note that 
Assumption \ref{assu:botelho} guarantees that 
Assumption \ref{assu:lorenz}\emph{\ref{assu:lorenz.ter}} holds.
Also, from Assumption \ref{assu:botelho} and the definitions of $ \beta_k $ and $ \beta $, we obtain 
$ \beta_k \geq \beta $ for all $ k \geq 0 $.
Then, since $ \beta_i \geq \beta $ ($ 0 \leq i \leq k $), 
we get, for all $ k \geq 0 $,
\begin{align*}
p_k = \left( \prod_{i=0}^k\,(1 - \beta_i)\right)^{-1} \geq
\left( \prod_{i=0}^k\,(1 - \beta)\right)^{-1} = \dfrac{1}{(1-\beta)^{k+1}}.
\end{align*}
Hence\,---\,see \eqref{eq:def.ubs02}\,---\,, for all $ k \geq 1 $,
\begin{align*}
\Lambda_k \geq \underline \lambda \eta \sum_{j=0}^{k-1}p_j 
\geq \sum_{j=0}^{k-1} \dfrac{\underline\lambda \eta }{(1-\beta)^{j+1}}
\geq \dfrac{\underline\lambda \eta }{(1-\beta)^k},
\end{align*}
where the term in the right-hand side of the third inequality is exactly the term in the summation 
with $ j = k-1 $.
\eproo

\mgap

\bprop[\bf Optimality; the case $ \seq{\eta_k} $ constant] \lab{pro:optstr02}
For all $ k \geq 1 $,
\begin{align} \lab{eq:natura}
-B_H \dist(\overline y_k, Q) \leq 
\emph{Gap}(\overline y_k, H, Q)\leq \dfrac{(1 - \beta)^k}{2\underline\lambda \eta}\left( D_X^2 + 
\sum_{j=0}^{k-1}\,p_{j-1}\delta_j \right),
\end{align}
where $ B_H > 0 $ and $ \beta \in (0, 1) $ are as in  \eqref{eq:def.ch} and \eqref{eq:def.beta}, respectively.
\eprop
\bproo
Since the assumptions of Proposition \ref{pro:optstr} are met, we conclude that
\eqref{eq:doc} holds, which gives \eqref{eq:natura} as a direct consequence of the second inequality in
\eqref{eq:belchior} and \eqref{eq:doc}.
\eproo

\mgap

In the next corollary, we specialize the result of Proposition \ref{pro:optstr02} for the case where no inertial effects are introduced in Algorithm \ref{alg:main}, in which case it reduces essentially to Algorithm~3.1 in~\cite{SamYou25}. 

\mgap

\bcoro[\bf Optimality; the case $ \seq{\eta_k} $ constant and $ \alpha_k \equiv 0 $] \lab{cor:optstr02b}
Under the assumptions of \emph{Proposition \ref{pro:optstr02}}, suppose that $ \alpha_k \equiv 0 $, i.e., suppose we have no inertial effects in \emph{Algorithm~\ref{alg:main}}.
Then, for all $ k \geq 1 $,
\begin{align} \lab{eq:peterlax}
-B_H \dist(\overline y_k, Q) \leq 
\emph{Gap}(\overline y_k, H, Q)\leq (1 - \beta)^k \left(\dfrac{D_X^2}{2\underline\lambda \eta}\right),
\end{align}
which means that, in this case, the dual gap $ \emph{Gap}(\overline y_k, H, Q) $ is reduced from above at a 
linear rate.
\ecoro
\bproo
By \eqref{eq:def.pd}, in this case we have $ \delta_k \equiv 0$, which gives \eqref{eq:peterlax} from
\eqref{eq:natura}.
\eproo

\mgap

\brema \lab{rem:rockafellar}
In the strongly monotone setting, sublinear rates were established  in \emph{\cite[Theorem 4.6]{SamYou25}} both for diminishing and constant regularization parameter sequences $ \seq{\eta_k} $ \emph{(}see \emph{Cases 1 and 2}\emph{)}. 
In the case of constant regularization $ \eta_k \equiv \eta > 0 $, i.e., \emph{Case 2}, $ \eta > 0 $ and the iteration index $ k $ are required to satisfy $ \eta \coloneqq \frac{(p + 1)\ln(k)}{\gamma \mu_H k} $ and $ \frac{k}{\ln(k)} \geq 5(p+1)\frac{L_H}{\mu_H} $
for some $ p \geq 1 $ \emph{(}in the notation of the latter reference\emph{)}. Also, in \emph{\cite[Theorem 4.6]{SamYou25}}, 
\emph{Case 3}, linear convergence rates were proved assuming that the feasible set $ Q $ is $ \sigma $-weakly sharp of order 
$ \mathcal{M} = 1 $ and $ \eta > 0 $ is chosen to satisfy $ \eta \leq \frac{\sigma}{2\norm{H(x^*)}} $ 
\emph{(}recall that the notation $ \alpha $ is used in place of $ \sigma $ in \emph{\cite{SamYou25}}\emph{)}.

In contrast, we established the linear rate \eqref{eq:peterlax} for the non-inertial version of \emph{Algorithm \ref{alg:main}}
\emph{(}which essentially reduces to \emph{\cite[Algorithm 3.1]{SamYou25}}\emph{)} without relying on any of those more restrictive assumptions imposed in \emph{\cite{SamYou25}}. 
This highlights that our analysis of \emph{Algorithm~1} remains relevant and of independent interest, even in the non-inertial setting, serving as a valuable complement to the results of~\emph{\cite{SamYou25}}.
\erema

\mgap

As opposite to Corollary \ref{cor:optstr02b}, in Corollaries \ref{cor:optstr02} and \ref{cor:optstr03} below, we present results on ``optimality'' for
Algorithm~\ref{alg:main} in the general (inertial) case. Before that, we state and prove the following lemma:

\blemm \lab{lem:kaizen}
For all $ k \geq 1 $,
\begin{align} \lab{eq:kaizen} 
\sum_{j=0}^{k-1}\,p_{j-1}\delta_j \leq 2kD_X^2.
\end{align}
\elemm
\bproo
First note that since the sequence $ \seq{p_{k-1}\alpha_k} $ is nonincreasing\,---\,see \eqref{eq:pkapk}\,---\,
and $p_{-1}\alpha_0 \leq 1 $ (because $ p_{-1} = 1$ and $ \alpha_0 \in [0,1 ] $), we conclude that
$ p_{k-1}\alpha_k \leq 1 $ for all $ k \geq 0 $. Moreover, 
Assumption \ref{assu:lorenz} clearly implies that
$ 0\leq \alpha_k \leq 1 $, and so $ \alpha_k^2 \leq \alpha_k $, for all $ k \geq 0 $. We also recall that, since
$ x_k, x_{k-1} \in X $, we have $ \norm{x_k - x_{k-1}}\leq D_X $ for all $ k \geq 0 $. As a consequence,
for $ k \geq 1 $ and $ 0 \leq j \leq k - 1 $, we derive that
\begin{align} \lab{eq:natura02}
p_{j-1}\delta_j = p_{j-1}\alpha_j(1+\alpha_j)\normq{x_j - x_{j-1}} \leq 2 \underbrace{p_{j-1}\alpha_j}_{\leq 1} D_X^2,
\end{align}
where the latter identity comes from the definition of $ \delta_j $. By summing up the inequalities in 
\eqref{eq:natura02} for $ j = 0, \dots, k-1 $, we obtain \eqref{eq:kaizen}. 
\eproo

\mgap

\bcoro[\bf Optimality; the case $ \seq{\eta_k} $ constant] \lab{cor:optstr02}
Under the assumptions of \emph{Proposition \ref{pro:optstr02}}, we have
\begin{align} \lab{eq:naturab}
-B_H \dist(\overline y_k, Q) \leq 
\emph{Gap}(\overline y_k, H, Q)\leq (k+1)(1 - \beta)^k \left(\dfrac{D_X^2}{\underline \lambda \eta}\right) \qquad
\mbox{for all} \quad k\geq 1.
\end{align}
\ecoro
\bproo
The proof follows directly from Proposition \ref{pro:optstr02} and Lemma \ref{lem:kaizen}.
\eproo

\mgap

\brema \lab{rem:shor}
We mention that although the convergence rate given in \eqref{eq:naturab} is not linear in the strict sense, it still resembles linear convergence in its behavior. This type of rate also appears in the analysis of inertial methods studied in \emph{\cite{WanFadOch25}}.
\erema

\mgap

In the next corollary, we infer about the number of iterations required to reduce (from above) the dual gap function  
$ \mbox{Gap}(\cdot, H, Q) $ at $ \overline y_k $ to below a given tolerance $ \veps > 0 $:

\mgap

\bcoro[\bf Optimality; the case $ \seq{\eta_k} $ constant] \lab{cor:optstr03}
Consider the assumptions of \emph{Proposition \ref{pro:optstr02}}, let $ \veps>0 $ be a given tolerance and let
$ k \geq 1 $. If
\begin{align} \lab{eq:logdx}
 k \geq \left\lceil \left(\dfrac{1}{1 - \overline\lambda^{\,2} L^2} 
+ \dfrac{1}{2\underline\lambda \eta \mu}\right)\log\left(\dfrac{(k+1)D_X^2}{\underline\lambda \eta \veps}\right) \right\rceil
\end{align}
then
\begin{align} \lab{eq:optdx}
-B_H \dist(\overline y_k, Q) \leq 
\emph{Gap}(\overline y_k, H, Q) \leq \veps.
\end{align}
\ecoro
\bproo
Indeed, assume that \eqref{eq:logdx} holds and \eqref{eq:optdx} does not hold, i.e., $ \text{Gap}(\overline y_k, H, Q) > \veps $. Then, in view of
\eqref{eq:naturab}, we obtain
\begin{align*}
(k+1)(1 - \beta)^k \left(\dfrac{D_X^2}{\underline \lambda \eta}\right) > \veps,
\end{align*}
which, in turn, leads to
\begin{align*}
k \log\left(\dfrac{1}{1-\beta}\right) < \log\left(\dfrac{(k+1)D_X^2}{\underline\lambda \eta \veps}\right).
\end{align*}
Since $ \log(\frac{1}{1-\beta}) \geq \beta $, we then conclude that $ k < \frac{1}{\beta} \log\left(\dfrac{(k+1)D_X^2}{\underline\lambda \eta \veps}\right) $ and, consequently,
\begin{align*}
k + 1 \leq \left\lceil \dfrac{1}{\beta}\log\left(\dfrac{(k+1)D_X^2}{\underline\lambda \eta \veps}\right) \right\rceil,
\end{align*}
which, due to the definition of $ \beta $ as in \eqref{eq:def.beta}, clearly contradicts \eqref{eq:logdx}.
This finishes the proof of the corollary.
\eproo

\mgap

\bprop[\bf Feasibility; the case $ \seq{\eta_k} $ constant] \lab{pro:feasstr02}
For all $ k \geq 1 $,
\begin{align} \lab{eq:gilmour}
0 
	\leq \emph{Gap}(\overline y_k, F, X) 
	\leq \dfrac{(1 - \beta)^k}{2\underline\lambda}\left( D_X^2 + \sum_{j=0}^{k-1}\,p_{j-1}\delta_j \right)
+ \eta\left(\dfrac{\overline\lambda C_H D_X}{\underline\lambda}\right),
\end{align}
where $ C_H $ and $ \beta \in (0, 1) $ are as in \eqref{eq:def.ch} and \eqref{eq:def.beta}, respectively.
\eprop
\bproo
Let $ k \geq 1 $. From Proposition \ref{pro:feasstr}, we have
\begin{align} \lab{eq:waters}
0 \leq \text{Gap}(\overline y_k, F, X) 
	& \leq \dfrac{\eta}{2\Lambda_k}\left( D_X^2 + \sum_{j=0}^{k-1}\,p_{j-1}\delta_j \right)
+ \dfrac{\eta^2 \overline\lambda}{\Lambda_k}\left(\sum_{j=0}^{k-1}\,p_j\right) C_H D_X.
\end{align}
From Lemma \ref{lem:raman},
\begin{align*}
 \dfrac{\eta}{2\Lambda_k} \leq \dfrac{(1- \beta)^k}{2\underline\lambda} \quad \mbox{and}
 \quad \sum_{j=0}^{k-1}\,p_j \leq \dfrac{\Lambda_k}{\underline\lambda \eta},
\end{align*}
which combined with \eqref{eq:waters} clearly yields \eqref{eq:gilmour}.
\eproo

\mgap

\bcoro[\bf Feasibility; the case $ \seq{\eta_k} $ constant] \lab{cor:feasstr02}
Under the assumptions of \emph{Proposition~\ref{pro:feasstr02}},
for all $ k \geq 1 $, we have
\begin{align} \lab{eq:james}
0 
	\leq \emph{Gap}(\overline y_k, F, X) 
	\leq (k + 1) (1 - \beta)^k \left( \dfrac{D_X^2}{\underline\lambda} \right)  + \eta\left(\dfrac{\overline\lambda C_H D_X}{\underline\lambda}\right).
\end{align}
\ecoro
\bproo
The proof follows directly from Proposition \ref{pro:feasstr02} and Lemma \ref{lem:kaizen}.
\eproo

\mgap

\bcoro[\bf Feasibility; the case $ \seq{\eta_k} $ constant] \lab{cor:feasstr03}
Under the assumptions of \emph{Proposition \ref{pro:feasstr02}}, assume additionally that
$ \eta > 0$ and $ \mathcal{D}_0 > 0 $ are as in \eqref{eq:etaconsb} and 
\eqref{eq:ubdz}, respectively, 
where $ \veps > 0 $ is a given tolerance and let $ k \geq 1 $.
If
\begin{align} \lab{eq:niver}
 k \geq \left\lceil \left(\dfrac{1}{1 - \overline\lambda^{\,2} L^2} 
+ \dfrac{\mathcal{D}_0}{\underline\lambda \mu \, \veps}\right)\log\left(\dfrac{2(k+1)D_X^2}
{\underline\lambda \, \veps}\right) \right\rceil
\end{align}
then
\begin{align} \lab{eq:shor}
0 \leq \emph{Gap}(\overline y_k, F, X) \leq \veps.
\end{align}
\ecoro

\mgap

\btheo[\bf Optimality and Feasibility; the case $ \seq{\eta_k} $ constant] \lab{the:ejinho}
Suppose the assumptions of \emph{Corollary \ref{cor:feasstr03}} hold and assume, additionally, that
$ \mathcal{D}_0 > \veps $. 
Then, for all
\begin{align} \lab{eq:niver02}
 k \geq \left\lceil \left(\dfrac{1}{1 - \overline\lambda^{\,2} L^2} 
+ \dfrac{\mathcal{D}_0}{\underline\lambda \mu \, \veps}\right)\log\left(\dfrac{2(k+1)D_X^2\mathcal{D}_0}
{\underline\lambda \, \veps^2}\right) \right\rceil,
\end{align}
the conclusions of \emph{Theorem~\ref{the:ofcons02}} hold: condition \eqref{eq:bertsekas} is satisfied, and, 
if $ Q $ is assumed to be $ \sigma $-weakly sharp of order $ \mathcal{M} \geq 1 $, then inequality \eqref{eq:bertsekas02} also holds.
\etheo
\bproo
Assume that \eqref{eq:niver02} holds. Then, using that $ \mathcal{D}_0/\veps > 1 $ and 
$ \eta = \veps/(2\mathcal{D}_0) $\,---\,see \eqref{eq:etaconsb}\,---\,,  we conclude that \eqref{eq:logdx}
and \eqref{eq:niver} also hold. So \eqref{eq:bertsekas} follows from Corollaries \ref{cor:optstr03} and 
\ref{cor:feasstr03}. To prove \eqref{eq:bertsekas02}, under the assumption that 
$ Q $ is $ \sigma $-weakly sharp of order $ \mathcal{M} \geq 1 $, use Proposition \ref{pro:optm02}
and~\eqref{eq:shor}.
\eproo

\mgap

We finish this section with two remarks about Theorem \ref{the:ejinho}:

\brema \lab{rem:suasuna}
\begin{enumerate}[label = \emph{(\roman*)}]
\item The (additional) assumption $ \mathcal{D}_0 > \veps $, which is not particularly restrictive and can be removed, allows for a shorter proof of \emph{Theorem~\ref{the:ejinho}}.
\item We observe that, although the bound on the right-hand side of \eqref{eq:niver02} depends on the iteration index 
$ k $, this dependence is dampened by the logarithmic function. Moreover, with respect to the precision
$ \veps > 0 $ and the diameter $ D_X > 0 $, the bound is, overall, dominated by
$ O\left(\frac{\mathcal{D}_0}{\veps}\right)$,
which reflects a sublinear rate of the form 
$ O\left(\frac{1}{k}\right)$. 
This stands in contrast to the (almost) linear convergence result established in \emph{Corollary~\ref{cor:optstr02}} 
\emph{(}see also \emph{Remark~\ref{rem:shor}}\emph{)}.
The emergence of this final sublinear-flavored rate—when both optimality and feasibility are considered—is essentially due to the presence of the constant term
$ \eta \overline\lambda C_H D_X/\underline\lambda $
in \eqref{eq:james}, which effectively “kills” the (almost) linear feasibility rate introduced by the term 
$ (k+1)(1 - \beta)^k $
in \emph{Corollary~\ref{cor:feasstr02}}.
Since the optimality result \emph{(}namely, the first inequality \emph{(}lower bound\emph{)} of \emph{Corollary~\ref{cor:optstr02}}\emph{)} also depends on feasibility estimates \emph{(}via \emph{Proposition~\ref{pro:optm02}}\emph{)}, the overall convergence rate becomes \emph{(}essentially\emph{)} sublinear, as stated in \emph{Theorem~\ref{the:ejinho}}.m
\end{enumerate}
\erema

\section{Numerical experiments} \lab{sec:ne}

In this section, we present some preliminary numerical experiments to show the effectiveness of the inertia. Our experiments were conducted in a Jupyter notebook executed on the Google Colab environment with Python 3.10.12. In all examples, we compare our algorithm IneIREG (Algorithm~\ref{alg:main}) with IR-EG$_{\verb|m,m|}$ proposed in~\cite{SamYou25}, that
will be simply called~IREG here, which corresponds to IneIREG with $\alpha_k = 0$ for all~$k$.

\subsection{Example 1: finding the best Nash equilibrium}

Following \cite[Subsection 6.1]{SamYou25}, we consider the simple two-player zero-sum game given by
\begin{align} \label{eq:game01}
\begin{cases}
\displaystyle{\min_{x_1}} \,f_1(x_1, x_2) \coloneqq 20 - 0.1 x_1 x_2 + x_1 \quad \mbox{s.t.} \quad x_1 \in X_1 \coloneqq 
\set{x_1 \mid 11 \leq x_1 \leq 60},\\[2mm]
\displaystyle{\min_{x_2}} \,f_2(x_1, x_2) \coloneqq -20 + 0.1 x_1 x_2 - x_1 \quad \mbox{s.t.} \quad x_2 \in X_2 \coloneqq
\set{x_2 \mid 10 \leq x_1 \leq 50}.
\end{cases}
\end{align}
The above problem is equivalent to VIP$(F, X)$ with $ F(x) \coloneqq Ax + b $, where
$ A \in \R^{2\times 2} $ and $ b \in \R^2 $ are given by
\begin{align}
A =
\begin{bmatrix}
0 & -0.1\\
0.1 & \phantom{-}0
\end{bmatrix}
\quad \mbox{and} \quad
b = 
\begin{bmatrix}
1 \\
0
\end{bmatrix},
\end{align} 
respectively, with $ X \coloneqq X_1 \times X_2 $.
For the upper-level VIP$(H, Q)$ we consider $ H \colon \R^2 \to \R^2 $ given by $ H(x) = x $ for all $ x \in \R^2 $.
Note that this corresponds to the problem of minimizing $ x \mapsto (1/2)\normq{x} $ over $ Q $ (recall that $ Q $ denotes the solution set of \eqref{eq:game01}). 
Note that $ L_F = 0.1 $ and $ L_H = 1 $ in this case.

In our experimental setup, both the IREG and IneIREG algorithms were implemented with a fixed stepsize set to $\lambda = 1$. For the regularization parameter $\eta_k$, we consider both the constant value of $\eta_k = 0.1$ and the diminishing value of $\eta_k = 0.1 / \sqrt{k+1} $. Furthermore, the IneIREG algorithm was tested using both constant inertial parameter $\alpha_k = 0.5$  and dynamic inertia as in~\eqref{eq:pen} with $m=1$, $\theta=0.1$, and $\rho=10^{-4}$. 

\begin{figure}[htbp]
    \centering
    \subfloat[Trajectories]{\includegraphics[width=0.45\textwidth]{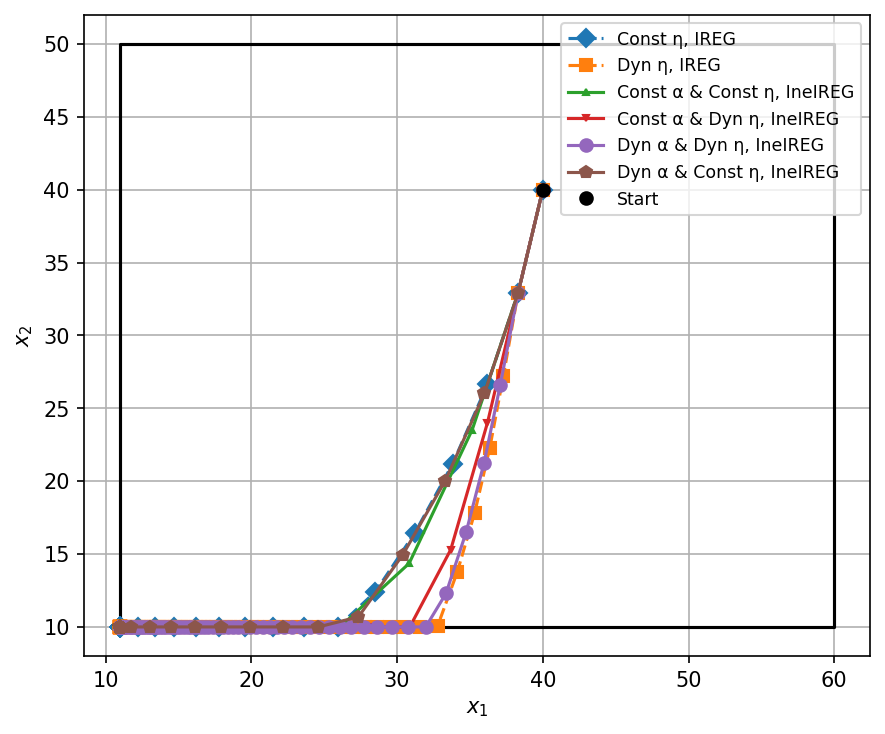}}
    \hfill
    \subfloat[Error vs. iteration]{\includegraphics[width=0.45\textwidth]{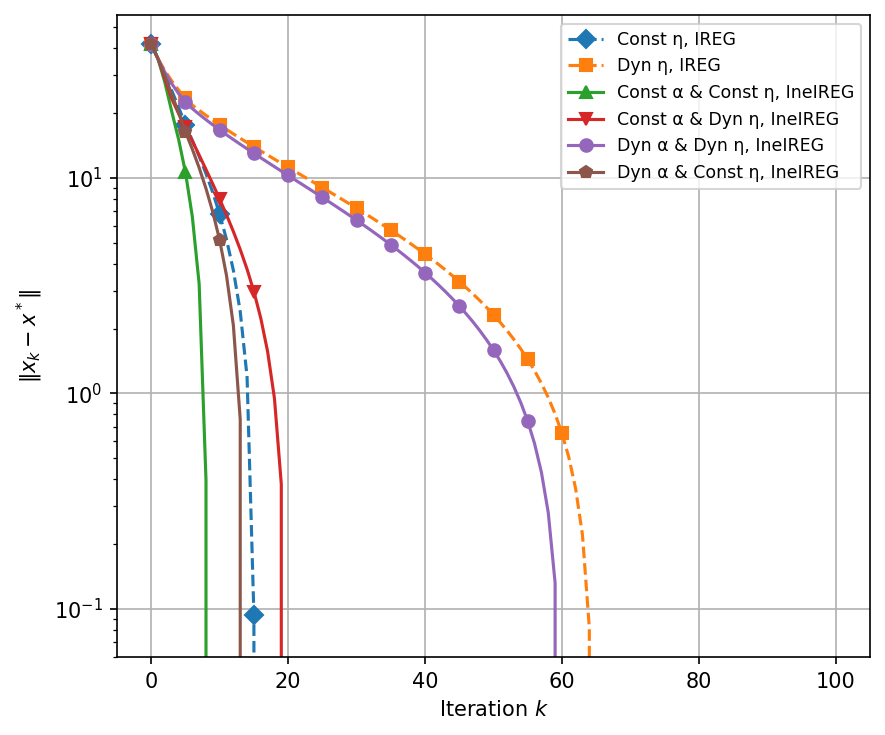}}
    \caption{Results for Example~1 using IREG and IneIREG}
    \label{fig: example1}
\end{figure}

The lines in Figure~\ref{fig: example1}(a) correspond to the trajectories of the iterates from the initial point $x_0 = [40, 40]^\top$ to the solution. Note that in all cases, we have convergence toward the same solution $x^* = [11, 10]^\top$. In Figure~\ref{fig: example1}(b), we plot how $\norm{x_k-x^*}$ changes through the iterations. We observe that the introduction of the inertial parameter $\alpha$ in IneIREG significantly accelerates convergence compared to the IREG algorithm. Comparing parameter settings, the combination of a constant $\eta$ and a constant $\alpha$ yields the best overall performance. These results confirm that inertia provides a substantial boost, with the specific parameterization being critical.

\subsection{Example 2: a simple random problem}

Now we consider another bilevel variational inequality problem given in~\cite[Example 2]{HieMou21}. 
The upper-level operator $H \colon \mathbb{R}^m \to \mathbb{R}^m$ is an affine mapping given by
\begin{align}
H(x) \coloneqq Qx + p.
\end{align}
% The vector $p \in \mathbb{R}^m$ has entries randomly drawn from a uniform distribution in $[-2, 2]$. 
The matrix $Q \in \mathbb{R}^{m \times m}$ is constructed to be symmetric and positive-definite by setting $Q = 0.5(Q_0 + Q_0^\top) + \tilde{m}I$, where $Q_0$ is a random matrix, $I$ is the identity matrix, and $\tilde{m}$ is a positive scalar (in our experiment, we just use $\tilde{m} = m$). This construction ensures that $H$ is strongly monotone and Lipschitz continuous.
 
The lower-level problem is a variational inequality defined by the feasible set $C$ and the operator $F \colon \mathbb{R}^m \to \mathbb{R}^m$, given by
\begin{align}
F(x) \coloneqq N^\top Nx + q,
\end{align}
where $N \in \mathbb{R}^{m \times m}$ and $q \in \mathbb{R}^m$.
% have entries randomly and uniformly drawn from the interval $[-2, 2]$. 
This choice guarantees that $F$ is monotone and Lipschitz continuous. The set $C \subset \mathbb{R}^m$ is a polyhedral convex set defined as
\begin{align}
C \coloneqq \{x \in \mathbb{R}^m_+ \mid Ex \leq f\}.
\end{align}
Here, the matrix $E \in \mathbb{R}^{10 \times m}$ contains entries randomly sampled from a uniform distribution in $[-1, 1]$, and the entries of the vector $f \in \mathbb{R}^{10}_+$ are all equal to~$m$ to ensure that $C$ is nonempty.

The experiments are conducted for dimension $m=100$. For all algorithms, we set the initial point as $x_0 = [1, 1, \ldots, 1]^\top$, which is feasible by construction. 
Moreover, $p$, $q$, $Q_0$, and $N$ have entries randomly chosen with a uniform distribution in $[-2, 2]$.
Since the exact solution to this problem is unknown, we check the performance of the algorithms by monitoring the convergence of the iterates via the metric $D_k := \|\bar{y}_{k+1} - \bar{y}_k\|$, where $\bar{y}_k$ is the ergodic mean given by~\eqref{eq:def.ubs}. We also use the metric
\begin{equation} \label{eq:infeasib}
  \phi(\bar{y}_{k}) := \| \max\{0, -\bar{y}_{k}\} \|^2 + \| \max\{0, -F(\bar{y}_{k})\} \|^2 + |\bar{y}_{k}^\top F(\bar{y}_{k})|
\end{equation}
to quantify the infeasibility. The common parameters of the methods are set as follows: the step-size is $\lambda = 1 / \|N^\top N\|$, the initial inertial parameter is $\alpha_0 = 0.1$, the parameters for the adaptive rule for $\alpha_k$ are $\theta = 0.99$ and $\rho = 10^{-8}$, and the initial regularization is given by $\eta_0 = 0.1$. 

\begin{figure}[htbp]
    \centering
    \subfloat[Suboptimality]{\includegraphics[width=0.45\textwidth]{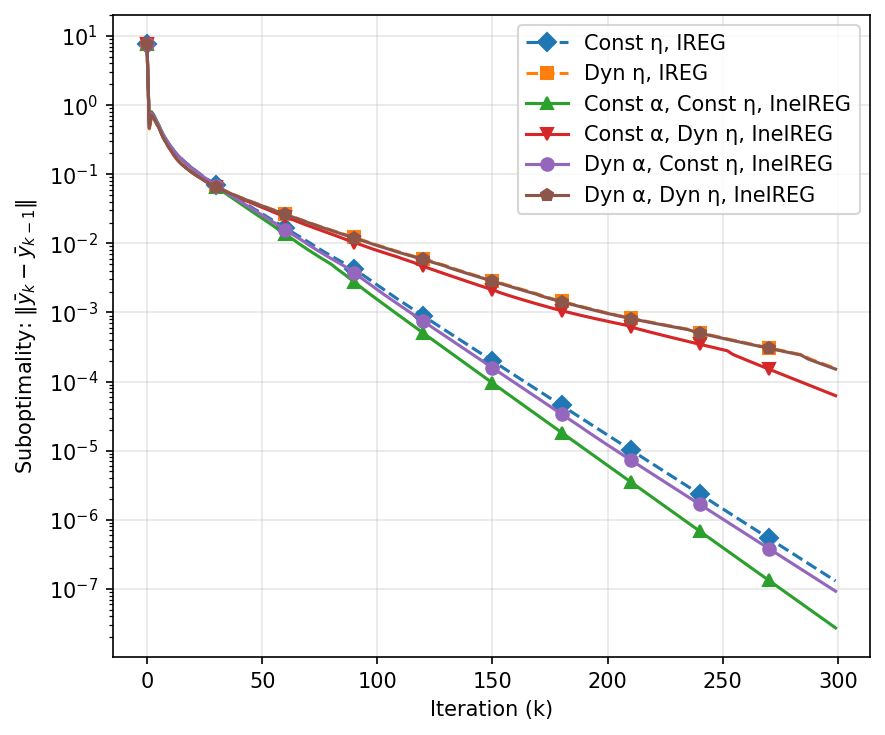}}
    \hfill
    \subfloat[Infeasibility]{\includegraphics[width=0.45\textwidth]{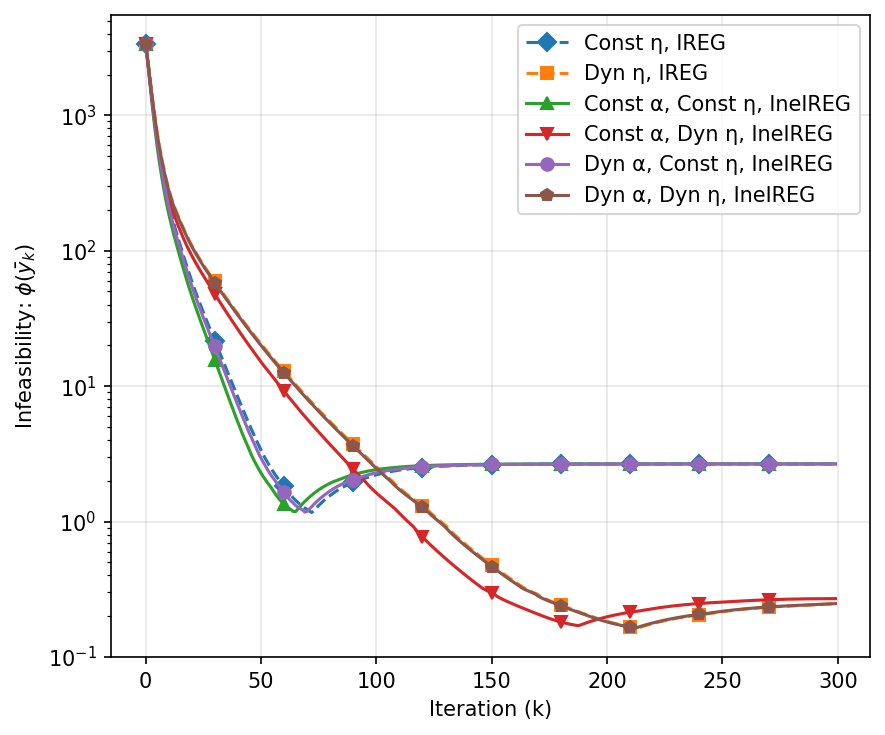}}
    \caption{Results for Example~2 using IREG and IneIREG}
    \label{fig: example2}
\end{figure}

Figures~\ref{fig: example2}(a) and \ref{fig: example2}(b) show the value of the metrics $D_k$ and $\phi(\bar{y}_{k})$, respectively, for each iteration~$k$. The numerical results indicate that the inertial parameter $\alpha$ leads to an improvement in performance with respect to the suboptimality metric, in particular when a constant $\alpha$ is employed. Furthermore, our experiments reveal that the most favorable outcomes are achieved when both the regularization parameter $\eta$ and the inertial parameter $\alpha$ are set to constant values. No significant difference in infeasibility was observed between IneIREG and IREG when applying the same regularization parameter settings.

\subsection{Example 3: a traffic equilibrium problem}

\begin{figure}[htbp]
    \centering 
    \includegraphics[width=0.7\textwidth]{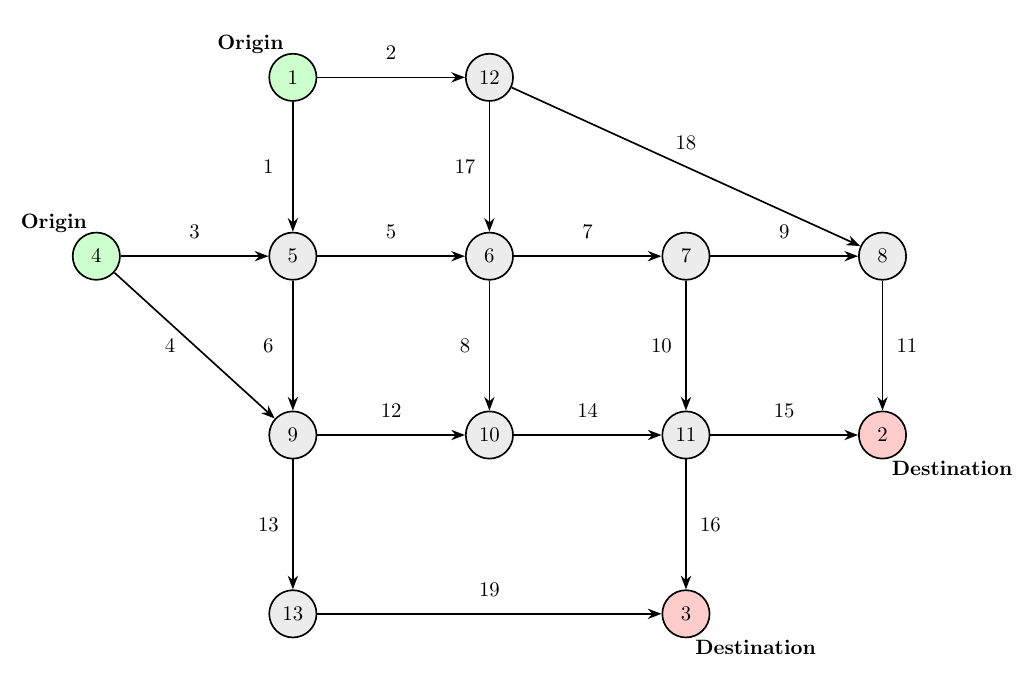} 
    \caption{Nguyen and Dupuis network} 
    \label{fig:network} 
\end{figure}

Similarly to~\cite[Subsection 6.2]{SamYou25}, we now consider the Nguyen and Dupuis traffic network equilibrium problem. The network, illustrated in Figure~\ref{fig:network}, consists of 13 vertices, 19 arcs, 25 paths, and 4 origin-destination (OD) pairs. Following the standard path-based formulation, the user equilibrium condition is modeled as a Nonlinear Complementarity Problem (NCP).
The goal is to find a state vector $x = [h; u]^\top \in \mathbb{R}^{29}$, where $h \in \mathbb{R}^{25}$ is the vector of traffic flows on the paths and $u \in \mathbb{R}^{4}$ is the vector of minimum travel costs between the OD pairs. The vector $x$ must satisfy the NCP conditions, formulated as:
\begin{align} \label{eq:ncp_traffic}
0 \le x \perp F(x) \ge 0,
\end{align}
where the operator $F \colon \mathbb{R}^{29} \to \mathbb{R}^{29}$ is defined by
\[
F(x) := \begin{bmatrix} C(h) - \Omega^\top u \\ \Omega h - d \end{bmatrix}.
\]
Here, $\Omega$ is the OD-path incidence matrix and $d$ is the given travel demand vector. The path cost vector $C(h)$ is derived from the arc travel times via the relation $C(h) = \Delta^\top c(\Delta h)$, where $\Delta$ is the arc-path incidence matrix. The travel time on each arc $a$, $c_a(\cdot)$, is given by the Bureau of Public Roads function:
\[
c_a(\mathcal{F}_a) := t_a^0 \left(1+ 0.15 \left(\frac{\mathcal{F}_a}{cap_a}\right)^{n_a}\right),
\]
where for each arc $a$, $\mathcal{F}_a$ is its flow, $t_a^0$ is the free-flow time, $cap_a$ is its capacity, and $n_a \ge 1$. It is shown in~\cite[Lemma~6.1]{SamYou25} that the mapping $F$ is monotone when $n_a \ge 0$ for all~$a$. 

The primary objective of this example is to find an equilibrium state that minimizes the total aggregated travel cost across the entire network. Define
\[
f(x) := \sum_{i=1}^{25} [C(h)]_i,
\]
which is convex when $n_a \ge 1$ for all arc~$a$ (see \cite[Lemma~6.2]{SamYou25}). The problem is thus stated as finding the best equilibrium, in the following sense:
\begin{align*}
\min \quad & f(x) \\
\text{s.t.} \quad & x \in \text{SOL}(\mathbb{R}^{29}_+, F), \nonumber
\end{align*}
where $\text{SOL}(\mathbb{R}^{29}_+, F)$ denotes the solution set of the NCP defined in~\eqref{eq:ncp_traffic}.

For all experiments, we take $n_a = 1$ and all the other problem's parameters are the same as~\cite[Section~6.2]{SamYou25}. We also take $x_0 = [1,\dots,1]^\top \in \mathbb{R}^{29}$ as the initial point. The common parameters for the methods are set as follows: the extragradient step-size is $ \lambda = 0.1 $, the initial inertial parameter is $\alpha_0 = 0.5$, and the initial regularization parameter is $\eta_0 = 0.1$. For cases with variable parameters, the decay exponent for $\eta_k$ is set as $b=0.5$, and the parameters for the adaptive rule for $\alpha_k$ are $\theta = 0.99$ and $\rho = 10^{-8}$. 
% We compare six different configurations by varying the strategies for the inertial parameter ($\alpha_k$) and the regularization parameter ($\eta_k$): (i) Fixed $\alpha$, Fixed $\eta$; (ii) Fixed $\alpha$, Variable $\eta$; (iii) Variable $\alpha$, Fixed $\eta$; (iv) Variable $\alpha$, Variable $\eta$; (v) No inertia ($\alpha_k=0$), Fixed $\eta$; and (vi) No inertia ($\alpha_k=0$), Variable $\eta$.

\begin{figure}[htbp]
    \centering
    \subfloat[Suboptimality]{\includegraphics[width=0.45\textwidth]{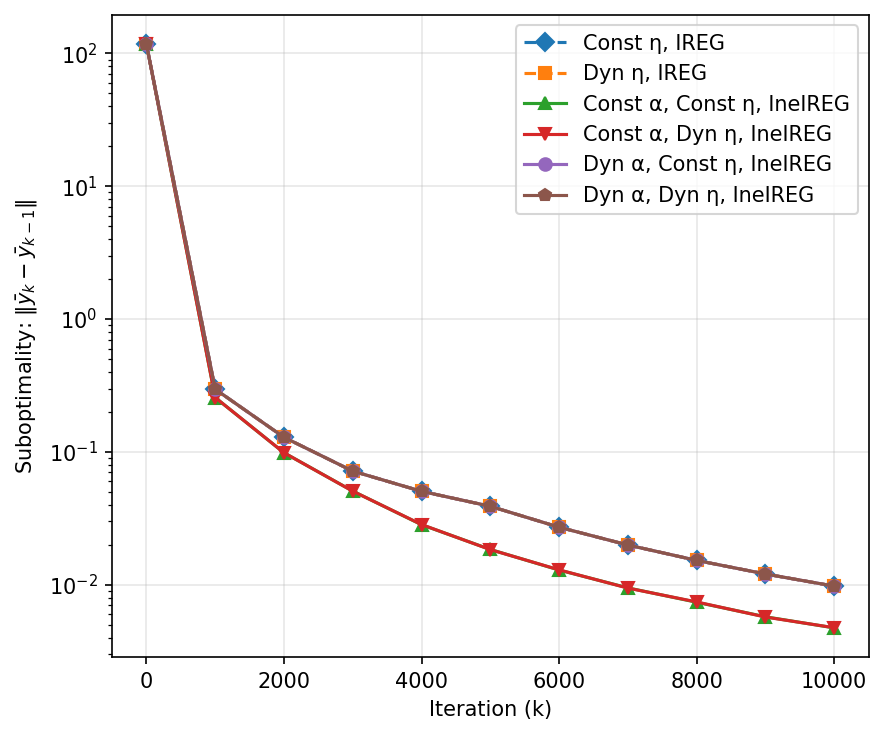}}
    \hfill
    \subfloat[Infeasibility]{\includegraphics[width=0.45\textwidth]{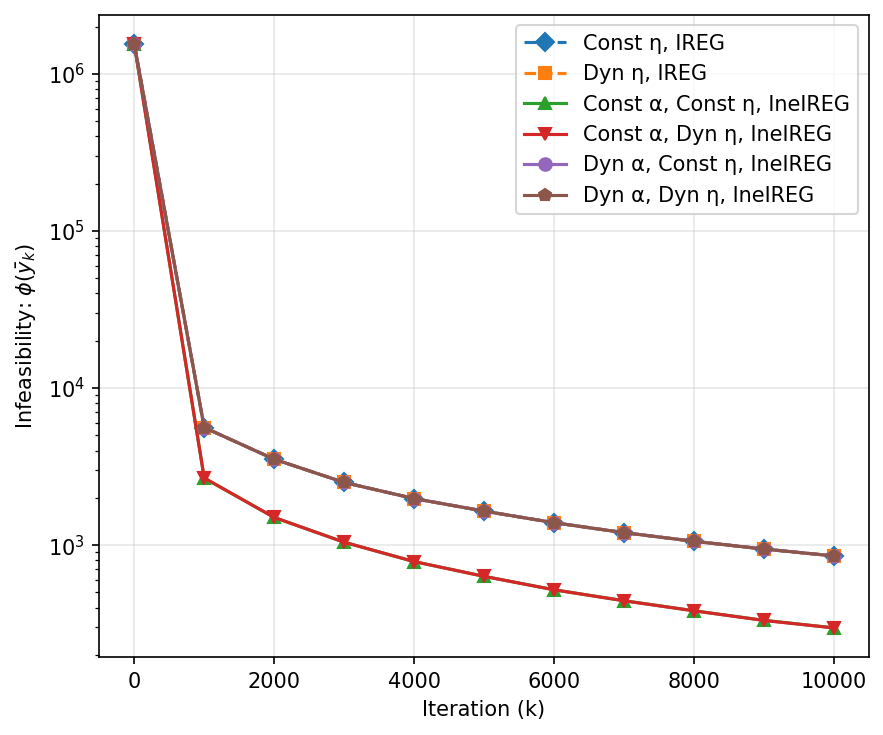}}
    \caption{Results for Example~3 using IREG and IneIREG}
    \label{fig: example3}
\end{figure}

We also illustrate the convergence behavior of the methods in terms of both suboptimality $\|\bar{y}_{k+1} - \bar{y}_k\|$ and infeasibility $\phi(\bar{y}_k)$ (given in~\eqref{eq:infeasib}).
The experimental results clearly show two performance groups based on the inertial parameter $\alpha$ as shown in Figure~\ref{fig: example3}.
Results using a constant $\alpha$ significantly outperform all others, demonstrating much faster convergence in both suboptimality and infeasibility. In contrast, the remaining methods (those with variable or no inertia) show slower convergence, underscoring the effectiveness of a simple, constant inertial strategy for this problem.

\section{Conclusions}
\label{sec:conclusions}

In this paper, we studied bilevel variational inequality problems, adding inertia in the iteratively regularized extragradient method, considered in~\cite{SamYou25}. We proved iteration-complexity bounds for monotone and strongly-monotone cases, refining also the results obtained in~\cite{SamYou25}. Our numerical experiments show that the addition of inertia may positively affect the performance of the algorithm. As a future work, we can analyze the particular case of VI-constrained optimization problems.

%\newpage
\appendix

{\small
\section*{Acknowledgments}
Part of this work was carried out while M.~M.~A.\ was visiting the Graduate School of Informatics at Kyoto University, Japan, whose hospitality is gratefully acknowledged.
}

\section{Auxiliary results}
 \lab{sec:app}

The proof of Theorem~\ref{the:projection} below can be seen, for instance, in~\cite[Proposition~2.2.1]{Ber03}. On the other hand, Lemma~\ref{lem:min} shows the minimizer of a strictly convex quadratic function under a linear inequality constraint.

\begin{theorem}[The projection theorem] \lab{the:projection}
Let $ X $ be a nonempty closed and convex subset of $ \HH $. Then, for all $ x \in \HH $, there exists a unique
element $ P_X(x) $ in $ X $, called the orthogonal projection of $ x $ onto $ X $, satisfying
$ \norm{x - P_X(x)}\leq \norm{x - y} $ for all $ y \in X $.
Moreover, $ z \in X $ is equal to $ P_X(x) $ if and only if
\[
\Inner{x - z}{x' - z}\leq 0 \quad \mbox{for all} \quad x' \in X.
\] 
\end{theorem}

\mgap
\mgap

\begin{lemma}\label{lem:min}
 Given $ a, b\geq 0 $ with $ a + b > 0 $ and $ c\in \mathbb{R}_+ $, we have
 \begin{align*}
  \min\{a s^2 + b t^2  \mid s,t\geq 0\;\;\mbox{and}\;\; s + t \geq c\} = \frac{ab}{a + b} c^2.
 \end{align*}
\end{lemma}

\section{Proof of Proposition \ref{pro:inv}} \lab{sec:app02}
\bproo
Let $ x \in X $ and $ k \geq 0 $. Using the identity $ \normq{a} - \normq{b} = \normq{a - b} + 2\inner{a - b}{b} $ with 
$ a \coloneq w_k - x $ and $ b \coloneq x_{k+1} - x $, we obtain
\begin{align} \lab{eq:ab}
\normq{w_k - x} - \normq{x_{k+1} - x} = \normq{w_k - x_{k+1}} + 2\inner{w_k - x_{k+1}}{x_{k+1} - x}.
\end{align}
Likewise, repeating the argument with $ a \coloneq w_k - y_k $ and $ b \coloneq x_{k+1} - y_k $, we also find
\begin{align} \lab{eq:ab2}
\normq{w_k - y_k} - \normq{x_{k+1} - y_k} = \normq{w_k - x_{k+1}} + 2\inner{w_k - x_{k+1}}{x_{k+1} - y_k}.
\end{align}
Subtraction of \eqref{eq:ab} and \eqref{eq:ab2} yields
\begin{align} \lab{eq:ab3}
\normq{w_k - x} - \normq{x_{k+1} - x} = \normq{w_k - y_k} - \normq{x_{k+1} - y_k} 
+ 2\inner{w_k - x_{k+1}}{y_k - x}.
\end{align}
Next we will find a lower bound to the term  $ 2\inner{w_k - x_{k+1}}{y_k - x} $ from \eqref{eq:ab3}.
To this end, first define $ G\colon \Dom F \cap \Dom H \to \HH $ by $ G(\cdot) = F(\cdot) + \eta_k H(\cdot) $
and note that $ G(\cdot) $ is monotone and Lipschitz continuous with (Lipschitz) constant
$ L_k > 0 $ as in \eqref{eq:def.Lk}. 
Using the definition of $ G(\cdot) $, the second identity in \eqref{eq:yx} and 
the projection theorem (see Theorem \ref{the:projection}), we obtain
\begin{align} \lab{eq:proj01}
\Inner{w_k - x_{k+1}}{x_{k+1} - x'} \geq \lambda_k\Inner{G(y_k)}{x_{k+1} - x'} 
\quad \text{for all} \quad x'\in X,
\end{align}
which, in turn, gives (after taking $ x' = x $ in \eqref{eq:proj01} and doing some simple algebra)
\begin{align} \lab{eq:inner01}
\nonumber
\inner{w_k - x_{k+1}}{y_k - x} & = \inner{w_k - x_{k+1}}{y_k - x_{k+1}} + \inner{w_k - x_{k+1}}{x_{k+1} - x}\\[2mm]
\nonumber
	& \geq \inner{w_k - x_{k+1}}{y_k - x_{k+1}} + \lambda_k\inner{G(y_k)}{x_{k+1} - x}\\[2mm]
	%\nonumber
	& = \inner{w_k - x_{k+1} - \lambda_k G(y_k)}{y_k - x_{k+1}} + \lambda_k\inner{G(y_k)}{y_k - x}.
\end{align}
Likewise, using the first identity in \eqref{eq:yx} combined with the projection theorem, we also find
\begin{align*}
\inner{w_k - \lambda_k G(w'_k) - y_k}{y_k - x'}\geq 0 \quad \text{for all} \quad x'\in X,
\end{align*}
which, in turn, gives (with $ x' = x_{k+1} $)
\begin{align} \lab{eq:inner02}
\nonumber
&\Inner{w_k - x_{k+1} - \lambda_k G(y_k)}{y_k - x_{k+1}} = \Inner{w_k - \lambda_k G(y_k) - y_k}{y_k - x_{k+1}} +\normq{y_k - x_{k+1}}\\[2mm]
\nonumber
& \hspace{1cm }= \Inner{w_k - \lambda_k G(w'_k) - y_k}{y_k - x_{k+1}}\\[2mm] 
\nonumber
& \hspace{2cm} + \lambda_k \Inner{G(w'_k) - G(y_k)}{y_k - x_{k+1}} +\normq{y_k - x_{k+1}}\\[2mm]
& \hspace{1cm} \geq \lambda_k \inner{G(w'_k) - G(y_k)}{y_k - x_{k+1}}
+\normq{y_k - x_{k+1}}.
\end{align}
In view of \eqref{eq:inner01} and \eqref{eq:inner02},
\begin{align} \lab{eq:inner03}
\nonumber
\inner{w_k - x_{k+1}}{y_k - x} & \geq 
\lambda_k \inner{G(w'_k) - G(y_k)}{y_k - x_{k+1}}\\[2mm]
& \hspace{3cm} +\normq{y_k - x_{k+1}} + \lambda_k\Inner{G(y_k)}{y_k - x}.
\end{align}
Direct use of \eqref{eq:inner03} into \eqref{eq:ab3} yields
\begin{align} \lab{eq:haar}
\nonumber
\normq{w_k - x} - \normq{x_{k+1} - x} & \geq \normq{w_k - y_k} + \normq{y_k - x_{k+1}} \\[2mm]
\nonumber
& \hspace{1cm}+ 2\lambda_k \Inner{G(w'_k) - G(y_k)}{y_k - x_{k+1}} + 2\lambda_k\Inner{G(y_k)}{y_k - x}\\[2mm]
\nonumber
& \geq \normq{w_k - y_k} - \lambda_k^2\normq{G(w'_k) - G(y_k)} + 2\lambda_k\inner{G(y_k)}{y_k - x}\\[2mm]
&\geq \Big( 1 - \lambda_k^2 L_k^2 \Big)\normq{w_k - y_k} + 2\lambda_k\inner{G(y_k)}{y_k - x},
\end{align}
where we also used the $L_k$-Lipschitz continuity of $ G(\cdot ) $ and the inequality
$ \norm{w'_k - y_k} \leq \norm{w_k - y_k }$ (recall that $ w'_k = P_\Omega(w_k) $ and $ y_k \in \Omega $).
To finish the proof of the proposition, use \eqref{eq:haar} and the definition of $ G(\cdot) $.
\eproo

%------------- BIBLIOGRAFIA --------------------------%

\bibliographystyle{plain}
%\bibliography{global,prox}

%\bibliography{C:/Users/maico/Dropbox/Bib_files/prox} % laptop
%\bibliography{/home/maicon/Dropbox/Bib_files/prox} % computador ufsc

\def\cprime{$'$}

\end{document}